\documentclass{amsart}
\usepackage{amssymb}

\newtheorem{theorem}{Theorem}[section]
\newtheorem{lemma}[theorem]{Lemma}
\newtheorem{corollary}[theorem]{Corollary}
\newtheorem{proposition}[theorem]{Proposition}

\theoremstyle{definition}
\newtheorem{definition}[theorem]{Definition}
\newtheorem{example}[theorem]{Example}

\newtheorem{criterion}[theorem]{Criterion}

\theoremstyle{remark}
\newtheorem{remark}[theorem]{Remark}
\newtheorem{notation}[theorem]{Notation}
\newtheorem*{claim}{Claim}
\newtheorem*{construc}{Construction}

\numberwithin{equation}{section}

\usepackage{pictex}

\begin{document}

\title{Topological mixing in $CAT\left(-1\right)$-spaces}
\author{Charalambos Charitos}
\address{Department of Mathematics, Agricultural University of Athens, 
75 Iera Odos,\linebreak Athens, Greece}
\email{bakis@auadec.aua.gr}
\thanks{This research was supported in part by Research Unit Grant 470}
\author{Georgios Tsapogas}
\address{Department of Mathematics, University of The Aegean, Karlovassi, 
Samos 83200, Greece}
\email{gtsap@aegean.gr}
\date{August 13, 1999 and, in revised form, May 18, 2000}

\begin{abstract}
If $X$ is a proper $CAT\left(  -1\right)  $-space and $\Gamma$ a
non-elementary discrete group of isometries acting properly discontinuously on
$X,$ it is shown that the geodesic flow on the quotient space $Y=X/\Gamma$ is
topologically mixing, provided that the generalized Busemann function has
zeros on the boundary $\partial X$ and the non-wandering set of the flow
equals the whole quotient space of geodesics $GY:=GX/\,\Gamma$ (the latter
being redundant when $Y$ is compact). Applications include the proof of
topological mixing for (A) compact negatively curved polyhedra, (B) compact
quotients of proper geodesically complete $CAT\left(  -1\right)  $-spaces by a 
one-ended group of isometries and (C) finite $n$-dimensional ideal polyhedra.
\end{abstract}

\maketitle

\section{Introduction and preliminary results}

The extensive study of the geodesic flow, primarily on Riemannian manifolds,
has been concerned, among other properties, with the establishment of
topological transitivity and topological mixing. For compact manifolds of
negative curvature, topological transitivity of the geodesic flow was proved
by Anosov in \cite{Ano}. Topological mixing, a stronger property, has been
shown for the class of compact manifolds with non-positive curvature by P.
Eberlein in \cite{Ebe2}. In this paper we establish topological mixing of the
geodesic flow in certain classes of spaces which are quotients of proper
$CAT\left(  -1\right)  $-spaces by a non-elementary discrete group of
isometries. A $CAT\left(  -1\right)  $-space is a geodesic metric space in
which every hyperbolic triangle is thinner than its associated comparison
triangle in the hyperbolic plane (for definitions and basic properties see
\cite{Bal}, \cite{B-G-H}, \cite{Bou} and \cite{Gro1}). The $CAT\left(
-1\right)  $ property is (locally) one among many possible generalizations to
singular spaces of the notion of negative curvature. Important examples of
$CAT\left(  -1\right)  $-spaces include Riemannian manifolds of sectional
curvature $\leq-1,$ metric trees and simply connected cell complexes of
negative curvature.

Throughout this paper we will use the letter $X$ to denote a proper
$CAT\left(  -1\right)  $-space. Let $\Gamma$ be a \textit{non-elementary
discrete} group of isometries of $X,$ i.e., the cardinality of the limit set of
the action of $\Gamma$ on $X$ is $>2$ (see section \ref{quo} below) and
consider the quotient metric space $Y=X$/ \negthinspace$\Gamma.$ Recall that
the space of geodesics $GX$ consists of all\ isometries $g:\mathbb{R}%
\rightarrow X$ and its topology is that of uniform convergence on
compact sets. The action of $\Gamma$ on $X$ induces an isometric action of
$\Gamma$ on the space $GX$ which is also properly discontinuous. Hence, the
space of geodesics $GY$ is defined as the quotient metric space $GX$/
\negthinspace$\Gamma$ (see section \ref{quo} below for details). We will use
the letter $p$ to denote both projections $X\rightarrow Y$ and $GX\rightarrow
GY.$ As the action of $\Gamma$ on $X$ is not necessarily free, observe that an
element $g\in GY$ is not a geodesic in the usual sense; it is just a
continuous map $g:\mathbb{R}\rightarrow Y$ for which there exists an isometry
$\overline{g}:\mathbb{R}\rightarrow X$ such that $g=p\circ\overline{g}.$ The
geodesic flow on $X$ is defined by the map
\[
\mathbb{R}\times GX\rightarrow GX
\]
where the action of $\mathbb{R}$ is given by a right translation, i.e., for all
$t\in\mathbb{R}$ and $g\in GX$, $\left(  t,g\right)  \rightarrow t\cdot g$
where $t\cdot g:\mathbb{R}\rightarrow X$ is the geodesic defined by
$\left(  t\cdot g\right)  \left(  s\right)  =g\left(  s+t\right)
,s\in\mathbb{R} $. If $t\in\mathbb{R}$ and $g\in GY$ define the geodesic flow
on $GY$ by setting
\[
t\cdot g=p\bigl(t\cdot \overline{g}\bigr)%
\]
where $\overline{g}$ is any lift of $g$ in $GX.$

\begin{definition}
\label{mix}The geodesic flow $\mathbb{R}\times GY\rightarrow GY$ is
topologically mixing if given any open sets $\mathcal{O}$ and $\mathcal{U}$ in
$GY$ there exists a real number $t_{0}>0$ such that for all $\left|  t\right|
\geq t_{0},$ $t\cdot \mathcal{O}\cap\mathcal{U}\neq\emptyset. $
\end{definition}

A point $g$ in $GY$ belongs to the \textit{non-wandering set} $\Omega$ of the
geodesic flow $\mathbb{R}\times GY\rightarrow GY$ if there exist sequences
$\left\{  g_{n}\right\}  \subset GY$ and $\left\{  t_{n}\right\}
\subset\mathbb{R},$ such that $t_{n}\rightarrow\infty,\ g_{n}\rightarrow g$
and $t_{n}\cdot g_{n}\rightarrow g.$

The \textit{generalized Busemann function} is a continuous function
$\alpha:\left(  \partial X\cup X\right)  \times X\times X\rightarrow
\mathbb{R}$ whose restriction on $X\times X\times X$ is given by
\[
\alpha\left(  y,x,x^{\prime}\right)  :=d\left(  x^{\prime},y\right)  -d\left(
x,y\right)
\]
for $\left(  \xi,x,x^{\prime}\right)  \in X\times X\times X$ (see section
\ref{BusSSS} for a precise definition). A main result of this paper, from
which we deduce topological mixing for the classes of spaces (A), (B) and (C)
explained below, is the following:

\begin{theorem}
\label{main}Let $X$ be a proper $CAT\left(  -1\right)  $-space and $\Gamma$ a
non-elementary discrete group of isometries of $X.$ Assume 

(1)
$\forall x,x^{\prime}\in X$ there exist $\xi\in\partial X$ such that
$\alpha\left(  \xi,x,x^{\prime}\right)  =0$, and 

(2) the
non-wandering set $\Omega$ equals $GY.$ \newline Then the geodesic flow on the
quotient space $Y=X$/ \negthinspace$\Gamma$ is topologically mixing.
\end{theorem}

\begin{remark}
If the quotient space $Y=X$/ \negthinspace$\Gamma$ is compact, the limit set
$\Lambda\left(  \Gamma\right)  $ of the action of $\Gamma$ on $X$ is the whole
boundary $\partial X$ (a proof of this is included in the proof of corollary
\ref{Bestvina} below). The property $\Lambda\left(  \Gamma\right)  =\partial
X$ is equivalent to $\Omega=GY$ (see proposition \ref{omega} below). Hence, if
$Y$ is compact, the second assumption in the above theorem, which is required
for the proof of topological mixing even in the case of manifolds with
non-negative curvature, is redundant.
\end{remark}

The above theorem implies that the geodesic flow is topologically mixing for
the following classes of spaces:

\begin{enumerate}
\item [(A)]compact negatively curved polyhedra, i.e., a finite union of
hyperbolic simplices glued together isometrically along faces of the same
dimension so that, with the induced metric, it has curvature $\leq-1.$

\item[(B)] compact quotients $X/\Gamma$ where $X$ is a proper
geodesically complete $CAT$ $\left(  -1\right)  $-space and $\Gamma$ a discrete
one-ended
group of isometries of $X.$

\item[(C)] $Y$ is an $n$-dimensional ideal polyhedron, i.e., a finite union of
ideal hyperbolic $n$-polytopes glued together isometrically along their
$\left(  n-1\right)  $-faces with at least two germs of polytopes along each
$\left(  n-1\right)  $-face so that with the induced metric $Y$ is a
complete length space of curvature $\leq-1.$
\end{enumerate}

The structure of this paper is as follows: in the present section we include
basic definitions and prove certain properties of $CAT\left(  -1\right)
$-spaces needed in the sequel. Moreover, the notion of a non-elementary group
$\Gamma$ of isometries of $X$ is explained and the action of such $\Gamma$ on
$GX$ is analyzed. Finally, the first section is concluded by 
(counter)examples, which are $1$-dimensional simplicial complexes, in which the
geodesic flow is not topologically mixing. These examples justify the fact
that $1$-dimensional simplices are excluded when proving topological mixing
for negatively curved polyhedra. In section \ref{BusSSS} Busemann functions
are discussed and we use them to study strong stable sets in the space of
geodesics. Although topological transitivity follows from topological mixing,
we show in section \ref{toptransitivity} that the geodesic flow on $Y $ is
topologically transitive because we need this property in the proof of theorem
\ref{main} given in section \ref{mainproof}. Finally, in section
\ref{application} the above-mentioned classes (A, B and C) of spaces are
defined and it is shown that all assumptions posited in theorem \ref{main}
above are satisfied.

\subsection{Preliminaries on $CAT\left(  -1\right)  $-spaces}\label{prelim}

The letter $X$ will always denote a proper $CAT\left(  -1\right)  $-space$.$
Recall that a metric space is proper if the closed balls are compact. For
definitions and basic properties of $CAT\left(  -1\right)  $-spaces we refer
the reader to \cite{Bal}, \cite{B-G-H} and \cite{Gro1}. We recall here basic
properties of the spaces $GX$ and $\partial X.$ $GX$ consists of
all\ isometric maps $g:\mathbb{R}\rightarrow X$ and its topology is the
topology of uniform convergence on compact sets. In addition, we note here
that $GX$ is metrizable and the metric is given by the formula
\begin{equation}
d_{GX}\left(  g_{1},g_{2}\right)  :=\int\nolimits_{-\infty}^{+\infty
}e^{-\left|  t\right|  }d\bigl(g_{1}\left(  t\right)  ,g_{2}\left(  t\right)
\bigr)\,dt  .\label{ggmetric}%
\end{equation}
If $g\in GX$ we will denote by $-g$ the geodesic defined by $\left(
-g\right)  \left(  s\right)  =g\left(  -s\right)  $ and, similarly, if
$\mathcal{A}\subset GX,$ then $-\mathcal{A}:=\left\{  -g\bigm\vert
g\in\mathcal{A}\right\}  .$

The (visual) boundary $\partial X$ of a $CAT\left(  -1\right)  $-space can be
defined, since $X$ is assumed to be proper, as the space of equivalence
classes of asymptotic geodesic rays starting at a fixed point in $X.$ If $g$
is a geodesic, we will denote by $g\left(  +\infty\right)  $ the boundary
point determined by the geodesic ray $g|_{\left[  0,+\infty\right)  }$ and
similarly for $g\left(  -\infty\right)  .$ We need the following two
conditions called (U) and (C), which are standard for studying flows (cf.
\cite{Kai}). Recall that two geodesic rays $g_{1},g_{2}$ (or geodesics) are
called \textit{asymptotic} if $d\bigl(g_{1}\left(  t\right)  ,g_{2}\left(
t\right)  \bigr)$ is bounded for all $t\in\mathbb{R}^{+}.$

\begin{description}
\item [$\mathrm{Condition\,\,(U)}$]For any two points $x_{1},x_{2}\in
X\cup\partial X$ there exists a unique geodesic joining them.

\item[$\mathrm{Condition\,\,(C)}$] For any two asymptotic geodesic rays (or
geodesics) $g_{1},g_{2}$ there exists a real number $d$ such that
\[
\mathrm{lim}_{t\rightarrow\infty}d\bigl(g_{1}\left(  t\right)  ,g_{2}\left(
t+c\right)  \bigr)=0.
\]
\end{description}

It is well known (see for example \cite{Ch}) that a $CAT\left(  -1\right)
$-space satisfies condition (U). This implies, in particular, that a complete
$CAT\left(  -1\right)  $-space is contractible. We next show that a proper,
complete $CAT\left(  -1\right)  $-space satisfies condition (C).

\begin{proposition}
\label{catx}A proper, $CAT\left(  -1\right)  $-space $X$ satisfies condition
(C). 
\end{proposition}

\begin{proof}Let $g_{1},g_{2}:\left[  0,\infty\right)  \rightarrow
X $ be two asymptotic geodesic rays. Denote by $\xi$ the common boundary point
$g_{1}\left(  +\infty\right)  =g_{2}\left(  +\infty\right)  .$ Let $\left\{
t_{n}\right\}  _{n\in\mathbb{N}}\subset\mathbb{R}$ be a sequence converging to
$+\infty.$ For each $n\in\mathbb{N}$, set $y_{n}=g_{2}\left(  t_{n}\right)  .$
The sequence $\left\{  s_{n}\right\}  _{n\in\mathbb{N}}\subset\mathbb{R}$
given by $s_{n}=d\bigl(g_{1}\left(  0\right)  ,y_{n}\bigr)$ converges to
$+\infty.$ Denote by $x_{n}$ the unique point on $\mathrm{Im}g_{1}$ such that
$d\bigl(g_{1}\left(  0\right)  ,x_{n}\bigr)=d\bigl(g_{1}\left(  0\right)
,y_{n}\bigr),$ i.e., $x_{n}=g_{2}\left(  s_{n}\right)  .$ For the reader's
convenience we have gathered all the above notation in figure 1.

\begin{figure}[tb]
\unitlength=1.00mm 
\begin{picture}(122.00,75.00)  \thicklines
\bezier{400}(10.00,11.00)(56.00,28.00)(119.00,30.00)
\bezier{400}(27.00,69.00)(41.00,52.00)(119.00,31.00)
\put(28.00,68.00){\circle*{2.00}}
\put(10.00,11.00){\circle*{2.00}}
\put(68.00,47.00){\circle*{2.00}}
\put(40.00,20.00){\circle*{2.00}}
\put(64.00,25.00){\circle*{2.00}}
\bezier{200}(27.00,69.00)(45.00,39.00)(64.00,25.00)%
\bezier{200}(27.00,69.00)(36.00,35.00)(40.00,20.00)%
\put(116.00,25.00){$\xi$} 
\put(21.00,71.00){$g_{1}(0)$}
\put(5.00,6.00){$g_{2}(0)$}
\put(39.00,16.00){$y_{n}=g_{2}(t_{n})$} 
\put(65.00,21.00){$y_m$}
\put(69.00,49.00){$x_{n}=g_{1}(s_{n})$}
\put(56.00,32.00){\circle*{2.00}} 
\put(57.00,34.00){$z_n$}
\label{fig2} 
\end{picture}
\caption{{}}%
\end{figure}

Let $c_{n}=t_{n}-s_{n},$ $n\in\mathbb{N}$. This sequence is increasing and
bounded above by $d\bigl(g_{1}\left(  0\right)  ,g_{2}\left(  0\right)  \bigr
).$ If $c$ is the real number such that $c_{n}\rightarrow c,$ we will show
that
\begin{equation}
\mathrm{lim}_{t\rightarrow\infty}d\bigl(g_{1}\left(  t\right)  ,g_{2}\left(
t+c\right)  \bigr)=0. \label{conducatx}%
\end{equation}
We will need the notion of the angle in $CAT\left(  -1\right)  $-spaces. We
refer the reader to \cite[Ch.I Sec.3]{Bal} for definitions and basic
properties. If $\left(  x,y,z\right)  $ is a geodesic triangle in $X,$ we
denote the angle subtended at $x$ by $\measuredangle_{x}\left(  y,z\right)  .$
Recall that if $\left(  \overline{x},\overline{y},\overline{z}\right)  $ is
the comparison triangle in the hyperbolic space $\mathbb{H}^{2}$ of the
geodesic triangle $\left(  x,y,z\right)  $, then
\[
\measuredangle_{x}\left(  y,z\right)  \leq\measuredangle\,_{\overline{x}%
}\left(  \overline{y},\overline{z}\right).
\]
We first show that
\begin{equation}
\measuredangle_{y_{n}}\bigl(g_{1}\left(  0\right)  ,g_{2}\left(  0\right)
\bigr)\rightarrow0\,\,\,as\,\,\,n\rightarrow\infty. \label{angle}%
\end{equation}
For each $n\in\mathbb{N}$, let $\left(  \overline{g_{1}\left(  0\right)
},\overline{g_{2}\left(  0\right)  },\overline{y_{n}}\right)  $ be the
comparison triangle in the hyperbolic space $\mathbb{H}^{2}$ of the geodesic
triangle $\left(  g_{1}\left(  0\right)  ,g_{2}\left(  0\right)
,y_{n}\right)  .$ Since $s_{n}\rightarrow\infty$ and $t_{n}\rightarrow\infty$
as $n\rightarrow\infty,$ it is not possible to have both angles
$\measuredangle\,_{\overline{g_{1}\left(  0\right)  }}\, \left( \overline
{g_{2}\left(  0\right)  },\overline{y_{n}}\right) $ and $\measuredangle
\,_{\overline{g_{2}\left(  0\right)  }}\, \left( \overline{g_{1}\left(  0\right)
},\overline{y_{n}}\right) $ converging to $0$ as $n\rightarrow\infty.$ Without
loss of generality, assume that $\measuredangle\,_{\overline{g_{2}\left(
0\right)  }}\, \left( \overline{g_{1}\left(  0\right)  },\overline{y_{n}}\right) $
is bounded away from zero for all $n.$ Then using the law of cosines
\[
\frac{\mathrm{sinh\,\,}d\left( \overline{g_{1}\left(  0\right)  }%
,\overline{g_{2}\left(  0\right)  }\right)}{\mathrm{sin\,\,}\left(
\measuredangle\,_{\overline{y_{n}}}\, \left( \overline{g_{1}\left(  0\right)
},\overline{g_{2}\left(  0\right)  }\right) \right)}=\frac{\mathrm{sinh\,\,}%
s_{n}}{\mathrm{sin\,\,}\left( \measuredangle\,_{\overline{g_{2}\left(
0\right)  }}\, \left( \overline{g_{1}\left(  0\right)  },\overline{y_{n}}\right
)\right)}%
\]
it follows that
\[
\measuredangle\,_{\overline{y_{n}}}\left(  \overline{g_{1}\left(  0\right)
},\overline{g_{2}\left(  0\right)  }\right)  \rightarrow0\,\,or,\,\,\,\pi
\,\,\,as\,\,\,n\rightarrow\infty.
\]
If $\measuredangle\,_{\overline{y_{n}}}\left(  \overline{g_{1}\left(
0\right)  },\overline{g_{2}\left(  0\right)  }\right)  \rightarrow\pi$, then
$s_{n}+t_{n}\rightarrow d\bigl(\overline{g_{1}\left(  0\right)  }%
,\overline{g_{2}\left(  0\right)  }\bigr)$ which is impossible (because
$\left\{  s_{n}\right\}  ,\left\{  t_{n}\right\}  \rightarrow+\infty$). Thus,
$\measuredangle\,_{\overline{y_{n}}}\left(  \overline{g_{1}\left(  0\right)
},\overline{g_{2}\left(  0\right)  }\right)  \rightarrow0.$ Since
$\measuredangle_{y_{n}}\left(  g_{1}\left(  0\right)  ,g_{2}\left(  0\right)
\right)  
\linebreak  
\leq\measuredangle\,_{\overline{y_{n}}}\left(
\overline{g_{1}\left(  0\right)  },\overline{g_{2}\left(  0\right)  }\right)
$, equation (\ref{angle}) is proved. \newline Our next step is to show that
\begin{equation}
d\left(  x_{n},y_{n}\right)  \rightarrow0\,\,\,as\,\,\,n\rightarrow
\infty. \label{xy}%
\end{equation}
The sequence of geodesic segments $\left[  g_{1}\left(  0\right)
,y_{n}\right]  $ converges to the geodesic ray $g_{1}$ uniformly on compact
sets. Thus for each $n\in\mathbb{N},$ we may find $m>n$ such that the
neighborhood (in the compact open topology) around $g_{1}$ determined by the
compact set $\left[  0,s_{n}\right]  $ and the positive number $1/n$ contains
the segment $\left[  g_{1}\left(  0\right)  ,y_{m}\right]  .$ In particular,
if $z_{n}$ is the unique point on $\left[  g_{1}\left(  0\right)
,y_{m}\right]  $ with $d\left(  g_{1}\left(  0\right)  ,z_{n}\right)  =s_{n}$
we have
\begin{equation}
d\left(  z_{n},x_{n}\right)  <1/n  .\label{zx}%
\end{equation}
In order to prove equation (\ref{xy}) above it suffices to show that
\begin{equation}
d\left(  z_{n},y_{n}\right)  \rightarrow0\,\,\,as\,\,\,n\rightarrow
\infty.\label{zy}%
\end{equation}
For each $n\in\mathbb{N}$, let $\left(  \overline{g_{1}\left(  0\right)
},\overline{y_{n}},\overline{y_{m}}\right)  $ be the comparison triangle of
the geodesic triangle $\left(  g_{1}\left(  0\right)  ,y_{n},y_{m}\right)  .$
Let $\overline{z_{n}}$ be the point corresponding to $z_{n}.$ Denote by
$\phi_{n}$ the angles $\measuredangle\,_{\overline{y_{n}}}\left(
\overline{g_{1}\left(  0\right)  },\overline{z_{n}}\right)  =\measuredangle
\,_{\overline{z_{n}}}\left(  \overline{g_{1}\left(  0\right)  },\overline
{y_{n}}\right)  .$ Apparently, $\phi_{n}<\pi/2$ for all $n.$ If $\left\{
\phi_{n}\right\}  ,$ or a subsequence, converges to $\phi,$ for some $\phi
<\pi/2,$ then using the facts
\begin{align*}
\measuredangle\,_{\overline{y_{n}}}\left(  \overline{g_{1}\left(  0\right)
},\overline{y_{m}}\right)   &  \rightarrow\pi\,\,\,\bigl(by\,\,\,(\ref{angle}%
)\bigr), \\
\measuredangle\,_{\overline{z_{n}}}\left(  \overline{g_{1}\left(  0\right)
},\overline{y_{m}}\right)   &  =\pi,
\end{align*}
it follows that
\[
\measuredangle_{\overline{y_{n}}}\left(  \overline{z_{n}},\overline{y_{m}%
}\right)  +\measuredangle_{\overline{z_{n}}}\left(  \overline{y_{n}}%
,\overline{y_{m}}\right)  >\pi/2+\pi/2,
\]
a contradiction. Therefore, $\phi_{n}\rightarrow\pi/2.$ Using this and the
second law of cosines we obtain that
\[
\mathrm{cosh}\,\,d\left(  \overline{z_{n}},\overline{y_{m}}\right)
\rightarrow1\,\,\,as\,\,\,n\rightarrow\infty;
\]
hence, $d\left(  \overline{z_{n}},\overline{y_{n}}\right)  \rightarrow0.$ By
comparison, $d\left(  z_{n},y_{m}\right)  \leq d\left(  \overline{z_{n}%
},\overline{y_{m}}\right)  $ which proves equation (\ref{zy}) and, in
consequence, proves equation (\ref{xy}).\newline We proceed now to show
equation (\ref{conducatx}). Since the function
\[
t\rightarrow d\bigl(g_{1}\left(  t\right)  ,g_{2}\left(  t+c\right)  \bigr)%
\]
is convex with respect to $t$ (see \cite[Ch. III]{B-G-H}), it suffices to show
that for each $\varepsilon>0$ there exists a positive real number $T=T\left(
\varepsilon\right)  $ such that
\[
d\bigl(g_{1}\left(  T\right)  ,g_{2}\left(  T+c\right)  \bigr)<\varepsilon.
\]
Let $\varepsilon>0$ be arbitrary. Choose $N\in\mathbb{N}$ such
\[
d\left(  x_{N},y_{N}\right)   <\varepsilon/2 ,\quad 
\left|  c_{N}-c\right|   <\varepsilon/2.\]
For the number $T=s_{N}$ we have
\[%
\begin{array}
[c]{lll}%
d\bigl(g_{1}\left(  T\right)  ,g_{2}\left(  T+c\right)  \bigr) &\!\!\! \leq &
\!\!\! d\bigl(g_{1}\left(  s_{N}\right)  ,g_{2}\left(  s_{N}+c_{N}\right)  \bigr
)+d\bigl(g_{2}\left(  s_{N}+c_{N}\right)  ,g_{2}\left(  s_{N}+c\right)  \bigr
)\\
 & \!\!\! = & \!\!\! d\bigl(g_{1}\left(  s_{N}\right)  ,g_{2}\left( 
t_{N}\right)  \bigr 
)+\left|  c_{N}-c\right| \\
 & \!\!\! = & \!\!\! d\left(  x_{N},y_{N}\right)  +\left|  c_{N}-c\right| 
<\varepsilon 
\end{array}
\]
which completes the proof of the proposition.\end{proof}

We will also need the following well-known lemma which asserts that the
projection of a point onto a geodesic always exists. For a proof see, for
example, \cite{C-T2}.

\begin{lemma}
\label{pro}Let $g$ be a geodesic in $G\widetilde{X}$ (or a geodesic segment)
and $x_{0}$ a point in $\widetilde{X}.$ There exists a unique real number $s$
such that $g\left(  s\right)  $ realizes the distance of $x_{0}$ from $\left.
\mathrm{Im\,}g,\right.  $ i.e., $dist\left(  x_{0},\mathrm{Im\,}g\right)
=d\left(  x_{0},g\left(  s\right)  \right)  .$
\end{lemma}

As usual, set $\partial^{2}X=\left\{  \left(  \xi,\eta\right)  \in\partial
X\times\partial X:\xi\neq\eta\right\}  .$ Condition (U) asserts that the fiber
bundle
\[
\rho:GX\rightarrow\partial^{2}X
\]
given by $\rho\left(  g\right)  =\bigl(g\left(  -\infty\right)  ,g\left(
+\infty\right)  \bigr)$ has a single copy of $\mathbb{R}$ as fiber. Moreover,
this bundle is trivial (see for example \cite[Th. 4.8]{Cham}). To define a
trivialization, let $x_{0}$ be a base point and let
\begin{equation}
H:G\widetilde{X}\overset{\approx}{\longrightarrow}\partial^{2}\widetilde
{X}\times\mathbb{R} \label{lization}%
\end{equation}
be the trivialization of $\rho$ with respect to $x_{0}$ defined by
\[
H\left(  g\right)  =\left(  g\left(  -\infty\right)  ,g\left(  +\infty\right)
,s\right)
\]
where $-s$ is the real number provided by lemma \ref{pro}.

It is shown in \cite[Prop.\ 4.8]{Cham} that the conjugation of the geodesic flow
with $H$ is simply the map
\begin{equation}
\left(  \xi_{1},\xi_{2},s\right)  \rightarrow\left(  \xi_{1},\xi
_{2},s+t\right)  ,\,for\,\,all\,\,\,\left(  \xi_{1},\xi_{2}\right)
\in\partial^{2}\widetilde{X}\,\,\,and\,\,\,s\in\mathbb{R}. \label{flow}%
\end{equation}

\subsection{The quotient space of geodesics}\label{quo}

In this section we define the space of geodesics for the quotient space $Y=X$/
\negthinspace$\Gamma$ and prove certain properties of it. We first recall the
notion of a non-elementary group of isometries. If $X$ is a $CAT\left(
-1\right)  $-space and $\Gamma$ a discrete group of isometries acting on $X,$
the \textit{limit set} $\Lambda\left(  \Gamma\right)  $ of the action of
$\Gamma$ is defined to be $\Lambda\left(  \Gamma\right)  =\overline{\Gamma
x}\cap\partial X,$ where $x$ is arbitrary in $X.$ The limit set has been
studied extensively (see \cite[Ch. II]{Coo1}, \cite[Ch. 2.1]{C-D-P} for a
detailed exposition) using the classification of the isometries of $X$ into
three types, namely, elliptic, parabolic and hyperbolic. If $\phi$ is
hyperbolic, then $\phi^{n}\left(  x\right)  $ converges to a point $\phi\left(
+\infty\right)  \in\partial X$ (resp. $\phi\left(  -\infty\right)  \in\partial
X $) as $n\rightarrow+\infty$ (resp. $n\rightarrow-\infty$) with $\phi\left(
+\infty\right)  \neq\phi\left(  -\infty\right)  $. Moreover,
\begin{equation}%
\begin{array}
[c]{l}%
\forall\,\,\xi\in\partial X\setminus\left\{  \phi\left(  +\infty\right)
\right\}  \bigl(resp.\,\,\partial X\setminus\left\{  \phi\left(
-\infty\right)  \right\}  \bigr)\Longrightarrow\\
\phi^{n}\left(  \xi\right)  \rightarrow\phi\left(  +\infty\right)  \,\,\bigl
(resp.\,\,\phi\left(  -\infty\right)  \bigr)\,\,as\,\,n\rightarrow
\infty\,\,\left(  resp.\,\,-\infty\right).
\end{array}
\label{7.2ofcoo}%
\end{equation}
The cardinality of the limit set is $0,1,2$ or infinite. A group $\Gamma$ acting
on a $CAT\left(  -1\right)  $-space $X$ is said to be \textit{non-elementary
}if the cardinality of $\Lambda\left(  \Gamma\right)  $ is infinite. In this
case, the following result is shown in \cite{Coo1}:
\begin{equation}
\left\{  \left(  \phi\left(  +\infty\right)  ,\phi\left(  -\infty\right)
\right)  :\phi\in\Gamma\,\,\,is\,\,\,hyperbolic\right\}
\,\,is\,\,dense\,\,in\,\,\Lambda\left(  \Gamma\right)  \times\Lambda\left(
\Gamma\right)  .\label{5.1ofcoo}%
\end{equation}
Note here that, as $X$ is assumed to be proper, discreteness of the
group $\Gamma$ is equivalent to requiring that $\Gamma$ acts properly
discontinuously on $X,$ i.e., for any compact$\,\,K\subset X$ the set $\left\{
\gamma\in\Gamma\bigm\vert\gamma K\cap K\neq\emptyset\right\}  $ is finite 
(see \cite[Th. 5.3.5]{Rat}). The following proposition is a well-known fact. We
include its proof here since we cannot find a reference for it.

\begin{proposition}
\label{discrete}Let $\Gamma$ be a group of isometries of $X$ acting properly
discontinuously on $X.$ Then $\Gamma$ acts by isometries and properly
discontinuously on the space of geodesics $GX.$
\end{proposition}

\begin{proof} We have assumed that $\Gamma$ acts by isometries on
$X.$ Therefore, if $f,g\in GX$ and $\gamma\in\Gamma$, we have
\[
\int_{-\infty}^{+\infty}e^{-\left|  t\right|  }d\left(  \gamma f\left(
t\right)  ,\gamma g\left(  t\right)  \right)  \,dt=\int_{-\infty}^{+\infty
}e^{-\left|  t\right|  }d\left(  f\left(  t\right)  ,g\left(  t\right)
\right)  \,dt
\]
which implies that $d_{GX}\left(  \gamma f,\gamma g\right)  =d_{GX}\left(
f,g\right)  .$ This shows that $\Gamma$ acts on $GX$ by isometries. Moreover,
we have assumed that $\Gamma$ acts properly discontinuously on $X,$ i.e.,
\begin{equation}
\forall\,\,compact\,\,K\subset X,\left\{  \gamma\in\Gamma\bigm\vert\gamma
K\cap K\neq\emptyset\right\}  \,\,is\,\,finite .\label{finite}%
\end{equation}
We proceed to show that $\Gamma$ acts properly discontinuously on $GX.$ Let
$\mathcal{K}$ be an arbitrary compact set in $GX.$ Set $d=diam\left(
\mathcal{K}\right)  $ and choose $g\in\mathcal{K}$ arbitrary. Using the
triangle inequality in $X$ one can show that
\[
d\bigl(\gamma g\left(  0\right)  ,g\left(  0\right)  \bigr)-2\left|  t\right|
\leq d\bigl(\gamma g\left(  t\right)  ,g\left(  t\right)  \bigr)\leq d\bigl
(\gamma g\left(  0\right)  ,g\left(  0\right)  \bigr)+2\left|  t\right|
\]
from which it follows, after integration, that
\begin{equation}
d\bigl(\gamma g\left(  0\right)  ,g\left(  0\right)  \bigr)-4\leq
d_{GX}\left(  \gamma g,g\right)  \leq d\bigl(\gamma g\left(  0\right)
,g\left(  0\right)  \bigr)+4 .\label{me}%
\end{equation}
As the space $X$ is assumed to be proper, the closure of the ball $B=B\bigl
(g\left(  0\right)  ,2d+4\bigr)$ centered at $g\left(  0\right)  $ and radius
$2d+4$ is compact and, by (\ref{finite}), the set
\[
A=\left\{  \gamma\in\Gamma\bigm\vert\gamma\overline{B}\cap\overline{B}%
\neq\emptyset\right\}  \,\,is\,\,finite.
\]
This together with equation (\ref{me}) implies that for all but a finite
number of elements $\gamma\in\Gamma,$%
\[
d_{GX}\left(  \gamma g,g\right)  \geq d\bigl(\gamma g\left(  0\right)
,g\left(  0\right)  \bigr)-4>2d+4-4=2d.
\]
Now let $f$ be an 
arbitrary element of $\mathcal{K}.$ Then, since $d_{GX}\left(
f,g\right)  <d=diam\left(  \mathcal{K}\right)  $, we have that for all but a
finite number of elements $\gamma\in\Gamma,$
\[%
\begin{array}
[c]{lll}%
d_{GX}\left(  \gamma f,g\right)  & \geq & \left|  d_{GX}\left(  \gamma
f,\gamma g\right)  -d_{GX}\left(  \gamma g,g\right)  \right| \\
& = & d_{GX}\left(  \gamma g,g\right)  -d_{GX}\left(  \gamma f,\gamma g\right)
\\
& > & 2d-d=d
\end{array}
\]
which implies that
\[
\forall\,f\in\mathcal{K}\Rightarrow\gamma f\notin\mathcal{K}%
\]
for all but a finite number of elements $\gamma\in\Gamma.$ In other words, the
set
\[
\left\{  \gamma\in\Gamma\bigm\vert\gamma\mathcal{K}\cap\mathcal{K}
\neq\emptyset\right\}
\]
is finite.\end{proof}

Define now $GY$ to be the orbit space $\left\{  \Gamma g\bigm\vert
g\in GX\right\}  $ of the action of $\Gamma$ on $GX.$ The space $GY$ can be
viewed as the set of all continuous functions $g:\mathbb{R}\rightarrow Y$ for
which there exists an isometry $\overline{g}:\mathbb{R}\rightarrow X$
satisfying $p\circ\overline{g}=g.$ By abuse of language, we will be calling
the elements of $GY$ geodesics in $Y.$

\begin{remark}
If, in addition, the action of $\Gamma$ on $X$ is free, so that $X$ would be
homeomorphic to the universal cover of $Y$ and the metric space $Y$ would have
curvature $\leq-1,$ then $GX$/ \negthinspace$\Gamma$ is the space of all local
geodesics $\mathbb{R}\longrightarrow Y$ (i.e., maps which are locally
isometric) which is, in fact, the natural definition for $GY.$
\end{remark}

Since the action of $\Gamma$ on $GX$ is properly discontinuous, each
$\Gamma$-orbit is a closed subset of $GX$ (cf. \cite[Th. 5.3.4]{Rat}). Using
this, the distance function
\[
d_{GY}:GY\times GY\rightarrow\mathbb{R}%
\]
defined by the formula
\[
d_{GY}\left(  \Gamma g,\Gamma f\right)  :=\mathrm{inf}\left\{  d\left(
x,y\right)  \bigm\vert x\in\Gamma g,y\in\Gamma f\right\}
\]
becomes a metric on $GY$ (cf. \cite[Th. 6.5.1]{Rat}). The topology induced by
the metric on $GY$ coincides with the quotient topology on $GY$ (see \cite[Th.
6.5.2]{Rat}). Moreover, it can be shown easily that the compact open topology
on $GY$ coincides with the quotient topology. Define the geodesic flow on $GY$
by the map
\[
\mathbb{R}\times GY\rightarrow GY:\left(  t,g\right)  \rightarrow t\cdot g
\]
where $t\cdot g=p\bigl(t \cdot \overline{g}\bigr)$ and $\overline{g}$ is
any lift of $g$ in $GX.$ It is easy to check that this definition does not
depend on the choice of the lift $\overline{g}.$

Since $\partial X$ is compact and $\partial^{2}X$ is an open subset of $\partial
X\times\partial X,$ $\partial^{2}X$ is separable. Moreover, $GX$ being, by
(\ref{lization}), homeomorphic to $\partial^{2}X\times\mathbb{R}, $ is also a
separable metric space. Thus, its continuous image $GY$ is separable, hence,
the metric space
\begin{equation}
GY\,\,is\,\,2{nd}\,\,countable .\label{secondcount}%
\end{equation}

\subsection{Non-mixing example}\label{cex}

We conclude this section by describing examples in which the geodesic flow is
not topologically mixing. In the example which is discussed in detail below
and also appears in \cite[Ch. II, Remark 3.6]{Bal}, the space $X$ is a
simplicial tree and the quotient space is a finite graph. In fact, as it was
pointed out to us by the referee, in the case
of a graph it is possible to characterize exactly when the geodesic flow is
mixing by looking at the lengths of the closed geodesics (see Remark
\ref{rfce} below).

\begin{example}
\label{ex}Let $Y$ be a plane graph homeomorphic to the figure-eight $S^{1}\vee
S^{1}$ consisting of seven vertices denoted $A,B,C,D,E,F,G$ and eight edges
$AB,BC,\! CD\! ,$ $DA,CE,EF,FG,GC$ all with length $1$ (see figure 2), let 
$X$ be its universal cover and $\Gamma$ the free group on two generators acting
on $X$ so that $Y=X$/ \negthinspace$\Gamma$. Then the geodesic flow on $GY$ is
not topologically mixing.
\end{example}

\begin{figure}[ptb]
\begin{picture}(310,140)
\thinlines
\put(60,63){\line(1,1){55}} 
\put(60,63){\line(1,-1){55}}
\put(170,63){\line(1,1){55}}
\put(170,63){\line(1,-1){55}}
\put(225,118){\line(1,-1){55}} 
\put(225,8){\line(1,1){55}} 
\put(115,118){\line(1,-1){55}}
\put(115,8){\line(1,1){55}} 
\put(60,63){\circle*{5}}%
\put(115,118){\circle*{5}} 
\put(170,63){\circle*{5}}%
\put(225,118){\circle*{5}} 
\put(115,8){\circle*{5}} 
\put(280,63){\circle*{5}}%
\put(225,8){\circle*{5}} 
\put(10,60){$A=g_1(0)$} 
\put(177,59){$C=g_2(0)$}
\put(125,59){$g_1(2)=$}
\put(119,119){$B=g_1(1)$} 
\put(119,0){$D=g_1(3)$} 
\put(285,60){$F=g_2(2)$} 
\put(229,0){$G=g_2(3)$}
\put(119,119){$B=g_1(1)$}
\put(229,119){$E=g_2(1)$} 
\label{}
\end{picture}
\caption{}
\end{figure}

\begin{proof}Observe first that as $X$ is a tree, it is a proper
geodesically complete $CAT\left(  -1\right)  $-space and its boundary is
totally connected. As $Y$ is compact, the second assumption of theorem
\ref{main} is also satisfied. Let $g_{1}:\mathbb{R}\rightarrow Y$ be the
closed geodesic with period $\omega=4$ satisfying $g_{1}\left(  0\right)  =A,$
$g_{1}\left(  1\right)  =B,$ $g_{1}\left(  2\right)  =C$ and $g_{1}\left(
3\right)  =D.$ Similarly, let $g_{2}:\mathbb{R}\rightarrow Y$ be the closed
geodesic with period $\omega=4$ satisfying $g_{2}\left(  0\right)  =C$ and
$g_{2}\left(  1\right)  =E.$ Observe that $g_{1},g_{2}$ are well defined by
the above requirements. Moreover, observe that for any $k\in\mathbb{Z}$
\[
g_{1}\left(  4k\right)  =A,g_{2}\left(  4k\right)  =C\,\,\,
\emph{and\,\,\,}%
g_{1}\left(  4k+1\right)  =B.
\]
Consider neighborhoods, in the compact open topology, $\mathcal{O}_{1}$ and
$\mathcal{O}_{2}$ of $g_{1}$ and $g_{2}$ respectively, determined by some
compact set in $\mathbb{R},$ say $\left[  -1/4,1/4\right]  ,$ and the number
$1/8$, i.e.,
\begin{equation}
h\in\mathcal{O}_{i}\Leftrightarrow d\bigl(h\left(  t\right)  ,g_{i}\left(
t\right)  \bigr)<1/8\,\,for\,\,all\,\,t\in\left[  -1/4,1/4\right]  ,\label{sex}%
\end{equation}
$i=1,2.$ We proceed to show that there exists a sequence $\left\{
t_{n}\right\}  \subset\mathbb{R}$ converging to $+\infty$ such that
$t_{n} \cdot \mathcal{O}_{1}\cap\mathcal{O}_{2}=\emptyset$ (cf. definition
\ref{mix} above). Let $\left\{  t_{n}\right\}  $ be the sequence $t_{n}=4n+1,$
$n\in\mathbb{N}.$

From (\ref{sex}) it is apparent that
\begin{equation}
h\in\mathcal{O}_{2}\Longrightarrow d\bigl(h\left(  0\right)  ,C\bigr)<1/8.
\label{o2}%
\end{equation}
Moreover, it is straightforward to check that if $h^{\prime}\in\mathcal{O}%
_{1}$, then $h^{\prime}\left(  t_{n}\right)  $ lies within a distance $1/8$ from
either $B$, $D$ or, $E$ or, $G.$ If $h\in t_{n} \cdot \mathcal{O}_{1} $,
then $h=t_{n} \cdot h^{\prime}$ for some $h^{\prime}\in\mathcal{O}_{1}, $
thus, $h\left(  0\right)  =h^{\prime}\left(  t_{n}\right)  .$ Hence,
\begin{equation}
h\in t_{n}\cdot \mathcal{O}_{1}\Longrightarrow dist\bigl(h\left(
0\right)  ,\left\{  B,D,E,G\right\}  \bigr)<1/8. \label{o1}%
\end{equation}
If $h\in t_{n}\cdot \mathcal{O}_{1}\cap\mathcal{O}_{2},$ combining
equations (\ref{o2}) and (\ref{o1}) above, we obtain that
\[%
\begin{array}
[c]{lll}%
1 & = & d\bigl(C,\left\{  B,D,E,G\right\}  \bigr)\\
& \leq &  d\bigl(h\left(  0\right)  ,C\bigr)+dist\bigl(h\left(  0\right)
,\left\{  B,D,E,G\right\}  \bigr)<1/8+1/8=1/4.
\end{array}
\]
This contradiction shows that $t_{n}\cdot \mathcal{O}_{1}\cap
\mathcal{O}_{2}=\emptyset$ for all $n\in\mathbb{N}$ completing the proof that
the geodesic flow on the figure-eight is not topologically mixing.\end{proof}

\begin{remark} \label{rfce}
A modification of the above argument can be used to  show that the geodesic
flow on any finite graph $Y$ (not homeomorphic to the circle $S^{1}$) is not
topologically mixing provided that the following condition holds:
\begin{equation}
\mbox{ for any two closed geodesics in }Y\mbox{ with periods }\ell_{1}\mbox{
and }\ell_{2},\mbox{ the ratio
}\frac{\ell_{1}}{\ell_{2}}\in\mathbb{Q}. \label{ratio}%
\end{equation}
Moreover, as it was pointed out to us by the referee, the following converse
statement can be shown: if $Y$ is any (finite or infinite) graph not
homeomorphic to the circle $S^{1},$ then the geodesic flow on $Y$ is
topologically mixing if
\begin{equation}
\mbox{ there exist two closed geodesics in }Y\mbox{ with periods }\ell
_{1}\mbox{ and }\ell_{2},\mbox{ such that
}\frac{\ell_{1}}{\ell_{2}}\notin\mathbb{Q}. \label{ifratio}%
\end{equation}
The proof of the latter statement utilizes the fact that $\frac{\ell_{1}}%
{\ell_{2}}\notin\mathbb{Q}$ implies that the set $\mathbb{Z}\ell
_{1}+\mathbb{Z}\ell_{2}$ is dense in $\mathbb{R},$ which, in turn, implies
that there exists a positive integer $M\in\mathbb{N}$ such that the distance
of the set $A=\left\{  n\ell_{1}+m\ell_{2}\bigm\vert-M\leq n,m\leq M\right\}
$ can be made arbitrarily small for all $s\in\left[  0,\ell_{1}+\ell
_{2}\right]  $ uniformly. Thus, the question of topological mixing on graphs
can be settled by looking at the subgroup of $\mathbb{R}$ generated by the
lengths of the closed loops.
\end{remark}

\section{Stable and strong stable sets}\label{BusSSS}

In this section we will define and study stable and strong stable sets in $GX
$ and $GY.$ For this we will use the \textit{generalized Busemann function}
(for more details see \cite[p.\ 27]{Bal}, \cite[Sec.\ 
2]{Kai}) whose definition we
recall briefly. As usual, $X$ will denote a proper $CAT\left(  -1\right)
$-space and $\Gamma$ a non-elementary group of isometries of $X$.

Define a function $\alpha:X\times X\times X\rightarrow\mathbb{R}$ by letting
\begin{equation}
\alpha\left(  \xi,x,x^{\prime}\right)  :=d\left(  x^{\prime},\xi\right)
-d\left(  x,\xi\right)  \label{generalized}%
\end{equation}
for $\left(  \xi,x,x^{\prime}\right)  \in X\times X\times X.$ It is shown in
\cite[Ch.\ II, Sec.\ 2]{Bal} that this function extends to a continuous function
\[
\left(  \partial X\cup X\right)  \times X\times X\rightarrow\mathbb{R}%
\]
denoted again by $\alpha,$ called the \textit{generalized Busemann function}.

This function, in fact, generalizes the classical Busemann function whose
definition makes sense in our context. To see this, let $\gamma:\left[
0,+\infty\right)  \rightarrow X$ be a geodesic ray; the Busemann function
associated to $\gamma$ is a function $b_{\gamma}$ on $X$ defined by
\[
b_{\gamma}\left(  x\right)  =\mathrm{lim}_{t\rightarrow\infty}\left[  d\left(
x,\gamma\left(  t\right)  \right)  -t\right].
\]
It is easy to see that for any $x\in X,$
\[%
\begin{array}
[c]{lll}%
\alpha\bigl(\gamma\left(  +\infty\right)  ,\gamma\left(  0\right)  ,x\bigr) &
= & \mathrm{lim}_{t\rightarrow\infty}\alpha\bigl(\gamma\left(  t\right)
,\gamma\left(  0\right)  ,x\bigr)\\
& = & \mathrm{lim}_{t\rightarrow\infty}\left[  d\left(  x,\gamma\left(
t\right)  \right)  -t\right]  =b_{\gamma}\left(  x\right)
\end{array}
\]
and, therefore, the Busemann function $b_{\gamma}$ coincides with $\alpha
\bigl(\gamma\left(  +\infty\right)  ,\gamma\left(  0\right)  ,\cdot\bigr),$
i.e., the restriction of $\alpha$ on $\left\{  \gamma\left(  +\infty\right)
\right\}  \times\left\{  \gamma\left(  0\right)  \right\}  \times X.$
Conversely, for arbitrary $\xi\in\partial X$ and $y\in X,$ the restriction
$\alpha\left(  \xi,y,\cdot\right)  \equiv\alpha|_{\left\{  \xi\right\}
\times\left\{  y\right\}  \times X}$ is simply the Busemann function
$b_{\gamma_{y\xi}}$ associated to the unique geodesic ray $\gamma_{y\xi}$ with
$\gamma_{y\xi}\left(  0\right)  =y$ and $\gamma_{y\xi}\left(  +\infty\right)
=\xi.$

The generalized Busemann function $\alpha$ is Lipschitz with respect to the
second and third variable with Lipschitz constant 1. The latter means, in
particular, that any Busemann function is Lipschitz with constant 1. To check
the Lipschitz property, let $\xi\in\partial X$ and choose a sequence $\left\{
z_{n}\right\}  \subset X$ such that $z_{n}\rightarrow\xi.$ Then for any fixed
$x\in X,$%
\[%
\begin{array}
[c]{lll}%
\left|  \alpha\left(  \xi,y,x\right)  -\alpha\left(  \xi,y^{\prime},x\right)
\right|  & = & \mathrm{lim}_{n\rightarrow\infty}\left|  d\left(
x,z_{n}\right)  -d\left(  y,z_{n}\right)  -d\left(  x,z_{n}\right)  +d\left(
y^{\prime},z_{n}\right)  \right| \\
& = & \mathrm{lim}_{n\rightarrow\infty}\left|  d\left(  y^{\prime}%
,z_{n}\right)  -d\left(  y,z_{n}\right)  \right|  \leq d\left(  y,y^{\prime
}\right).
\end{array}
\]
For $\xi\in X$ the calculation is analogous. Similarly, $\alpha$ can be shown
to be Lipschitz with respect to the third variable.

\begin{definition}
\label{ssset}We say that a geodesic $h\in GX$ belongs to the stable set
$W^{s}\left(  g\right)  $ of a geodesic $g$ if $g,h$ are asymptotic. Two
points $x,x^{\prime}\in X$ are said to be equidistant from a point $\xi
\in\partial X$ if $\alpha\left(  \xi,x,x^{\prime}\right)  =0.$

We say
that a geodesic $h\in GX$ belongs to the strong stable set $W^{ss}\left(
g\right)  $ of a geodesic $g$ if $h\in W^{s}\left(  g\right)  $ and $g\left(
0\right)  ,h\left(  0\right)  $ are equidistant from $g\left(  \infty\right)
=h\left(  \infty\right)  $.

Similarly, if $h,g\in GY,$ we say that
$h\in W^{ss}\left(  g\right)  $ $\bigl($respectively $W^{s}\left(  g\right)
\bigr)$ if there exist lifts $\overline{h},\overline{g}\in GX$ of $h,g$ such
that $\overline{h}\in W^{ss}\left(  \overline{g}\right)  $ $\bigl(%
$respectively $W^{s}\left(  \overline{g}\right)  \bigr).$
\end{definition}

The following proposition is a consequence of condition (C) and of the
properties of the $\alpha$ function.

\begin{proposition}
\label{asy}Let $f,g\in GX$ with $f\in W^{ss}\left(  g\right)  .$ Then
\[
\mathrm{lim}_{t\rightarrow\infty}d\bigl(f\left(  t\right)  ,g\left(  t\right)
\bigr)=0.
\]
\end{proposition}

\begin{proof} We first show that if $\alpha\bigl(\xi,f\left(
0\right)  ,g\left(  0\right)  \bigr)=0,$ where $\xi=f\left(  +\infty\right)
=g\left(  +\infty\right)  ,$ then
\begin{equation}
\alpha\bigl(\xi,f\left(  T\right)  ,g\left(  T\right)  \bigr)%
=0\,\,for\,\,all\,\,T\in\mathbb{R} .\label{t}%
\end{equation}
Fix $T>0$ (we work similarly for $T<0$). Choose a sequence $\left\{
t_{n}\right\}  \subset\mathbb{R}$, $t_{n}>T$ converging to $+\infty$. Then
\[%
\begin{array}
[c]{lll}%
\alpha\bigl(f\left(  t_{n}\right)  ,f\left(  0\right)  ,g\left(  0\right)
\bigr) & = & d\bigl(f\left(  t_{n}\right)  ,g\left(  0\right)  \bigr)-d\bigl
(f\left(  t_{n}\right)  ,f\left(  0\right)  \bigr)\\
& = & d\bigl(f\left(  t_{n}\right)  ,g\left(  0\right)  \bigr)-d\bigl(f\left(
t_{n}\right)  ,f\left(  T\right)  \bigr)-d\bigl(f\left(  T\right)  ,f\left(
0\right)  \bigr)\\
& = & \alpha\bigl(f\left(  t_{n}\right)  ,f\left(  T\right)  ,g\left(
0\right)  \bigr)-T.
\end{array}
\]
By taking the limits as $t_{n}\rightarrow\infty$ we have, by continuity of
$a,$ that
\[
\alpha\bigl(\xi,f\left(  T\right)  ,g\left(  0\right)  \bigr)=T.
\]
A similar calculation shows that
\[%
\begin{array}
[c]{lll}%
\alpha\bigl(g\left(  t_{n}\right)  ,f\left(  T\right)  ,g\left(  0\right)
\bigr) & = & \alpha\bigl(f\left(  t_{n}\right)  ,f\left(  T\right)  ,g\left(
T\right)  \bigr)+T
\end{array}
\]
which implies that $\alpha\bigl(\xi,f\left(  T\right)  ,g\left(  T\right)
\bigr)=\alpha\bigl(\xi,f\left(  T\right)  ,g\left(  0\right)  \bigr)-T=0.$
This completes the proof of equation (\ref{t}). A repetition of the argument
above asserts that for any $s\in\mathbb{R},$
\begin{equation}
\alpha\bigl(\xi,f\left(  T\right)  ,g\left(  T+s\right)  \bigr)%
=s\,\,for\,\,all\,\,T\in\mathbb{R} .\label{te}%
\end{equation}
Now let $c$ be the real number posited by Condition (C) making 
\[\mathrm{lim}_{t\rightarrow\infty}d\bigl(f\left(  t\right)  ,g\left(
t+c\right)  \bigr)=0.\]
 We show that $c=0$ concluding the proof of the proposition.
Assume on the contrary that $d\neq0.$ Let $T_{0}$ be large enough so that
\[
\left|  d\bigl(f\left(  T_{0}\right)  ,g\left(  T_{0}+c\right)  \bigr)\right|
<\left|  c\right|  /2.
\]
Choose a sequence $\left\{  t_{n}\right\}  \subset\mathbb{R}$ converging to
$+\infty$ with $t_{n}>T_{0}.$ Then
\[%
\begin{array}
[c]{lll}%
\left|  \alpha\bigl(f\left(  t_{n}\right)  ,f\left(  T_{0}\right)  ,g\left(
T_{0}+c\right)  \bigr)\right|  & = & \left|  d\bigl(f\left(  t_{n}\right)
,g\left(  T_{0}+c\right)  \bigr)-d\bigl(f\left(  t_{n}\right)  ,f\left(
T_{0}\right)  \bigr)\right| \\
& \leq & \left|  d\bigl(f\left(  T_{0}\right)  ,g\left(  T_{0}+c\right)
\bigr)\right|  <\left|  c\right|  /2
\end{array}
\]
and by taking the limit as $t_{n}\rightarrow\infty$ we have, by continuity of
$a,$ that 
\[\left|  \alpha\bigl(\xi,f\left(  T_{0}\right)  ,g\left(
T_{0}+c\right)  \bigr)\right|  \leq\left|  c\right|  /2,\]
a contradiction, by
equation (\ref{te}).\end{proof}

We will need the following
two lemmata concerning Busemann functions and strong stable sets.

\begin{lemma}
\label{buse}Let $\beta$ a geodesic of $X$ and $x\in X$ arbitrary. Then:

\textrm{(a)} The function $\alpha\bigl(\beta\left(  +\infty\right)
,x,\cdot\bigr):\mathrm{Im}\beta\rightarrow\mathbb{R}$ is an isometry.

\textrm{(b)} If $\gamma$ is any geodesic asymptotic with $\beta,$
then there exists a unique re-parametri\textit{\-}zation $\beta^{\prime}$ of
$\beta$ such that $\alpha\bigl(\beta\left(  +\infty\right)  ,\gamma\left(
0\right)  ,\beta^{\prime}\left(  0\right)  \bigr)=0,$ i.e. $\gamma\in
W^{ss}\left(  \beta^{\prime}\right)  .$ 

\textrm{(c)} Let $\gamma$ be
a geodesic ray in $X$ such that $\beta\left(  -\infty\right)  =\gamma\left(
+\infty\right)  .$ Then,
\[
\alpha\bigl(\gamma\left(  +\infty\right)  ,\gamma\left(  0\right)
,\beta\left(  t\right)  \bigr)=t+\alpha\bigl(\gamma\left(  +\infty\right)
,\gamma\left(  0\right)  ,\beta\left(  0\right)  \bigr).%
\]
In other words, the Busemann function $b_{\gamma}$ associated to $\gamma$ is
linear when restricted to $\mathrm{Im}\beta.$
\end{lemma}

\begin{proof}(a) Fix $t,t^{\prime}\in\mathbb{R}.$ Let $\left\{
x_{n}\right\}  \subset\mathrm{Im}\beta$ be a sequence converging to
$\beta\left(  +\infty\right)  .$ It is easily shown that for any $x\in X$ and
for all $n$ large enough $\bigl($namely, $\forall\,\,n$ for which
$x_{n}>\mathrm{max}\left\{  t,t^{\prime}\right\}  \bigr)$%
\[%
\begin{array}
[c]{lll}%
\left|  \alpha\left(  x_{n},x,\beta\left(  t\right)  \right)  -\alpha\left(
x_{n},x,\beta\left(  t^{\prime}\right)  \right)  \right|  & = & \left|
d\left(  x_{n},\beta\left(  t\right)  \right)  -d\left(  x_{n},\beta\left(
t^{\prime}\right)  \right)  \right| \\
& = & \left|  \beta\left(  t\right)  -\beta\left(  t^{\prime}\right)  \right|.
\end{array}
\]
Using the continuity of the $\alpha$ function and the fact that $x_{n}%
\rightarrow\beta\left(  +\infty\right)  $ we obtain that $\alpha\bigl(%
\beta\left(  +\infty\right)  ,x,\cdot\bigr)$ is an isometry on $\mathrm{Im}%
\beta.$

Part (b) follows from (a) by choosing $x=\gamma\left(
0\right)  $ and then defining $\beta^{\prime}\left(  t\right)  =\beta\left(
t+T\right)  $ where $T$ is the unique real number such that $\beta\left(
T\right)  $ is the inverse image of $0$ via the isometry $\alpha\bigl(%
\beta\left(  +\infty\right)  ,x,\cdot\bigr),$ i.e. $\alpha\bigl
(\beta\left(  +\infty\right)  ,x,\beta\left(  T\right)  \bigr)=0$.

(c)
Set $\xi=\beta\left(  -\infty\right)  =\gamma\left(  +\infty\right)  .$ Using
a sequence $\left\{  x_{n}\right\}  $ converging to $\xi$ and the continuity
of the $\alpha$ function it is easily shown that
\[
\alpha\bigl(\xi,\gamma\left(  0\right)  ,x\bigr)-\alpha\bigl(\xi,\beta\left(
0\right)  ,x\bigr)=\alpha\bigl(\xi,\gamma\left(  0\right)  ,\beta\left(
0\right)  \bigr),\,\,\forall\,\,x\in X.
\]
Hence, for arbitrary $t\in\mathbb{R}$ we have
\[
\alpha\bigl(\xi,\gamma\left(  0\right)  ,\beta\left(  t\right)  \bigr)%
=\alpha\bigl(\xi,\beta\left(  0\right)  ,\beta\left(  t\right)  \bigr)%
+\alpha\bigl(\xi,\gamma\left(  0\right)  ,\beta\left(  0\right)  \bigr).
\]
Pick $\left\{  t_{n}\right\}  \subset\mathbb{R},$ with $t_{n}\rightarrow
-\infty.$ Then,
\[%
\begin{array}
[c]{lll}%
\alpha\bigl(\xi,\beta\left(  0\right)  ,\beta\left(  t\right)  \bigr) & = &
\mathrm{lim}_{n\rightarrow\infty}\alpha\bigl(\beta\left(  t_{n}\right)
,\beta\left(  0\right)  ,\beta\left(  t\right)  \bigr)\\
& = & \mathrm{lim}_{n\rightarrow\infty}\bigl(t+\left|  t_{n}\right|  -\left|
t_{n}\right|  \bigr)=t.
\end{array}
\]
This completes the proof of the lemma.\end{proof}

\begin{lemma}
\label{contofsets} (a) For any $g\in GY$ and $c\in\mathbb{R},$ $\overline
{W^{ss}\left(  c\cdot g\right)  }=c\cdot \left(  \overline{W^{ss}\left(
g\right)  }\right)  .$

(b) Let $h_{1},g_{1}\in GY$ with $h_{1}\in
W^{ss}\left(  g_{1}\right)  $ and $\mathcal{O}_{1}\subset GY$ an open set
containing $h_{1}.$ Then there exists an open set $\mathcal{A}_{1}$ containing
$g_{1}$ such that for any $g\in\mathcal{A}_{1},$ $W^{ss}\left(  g\right)
\cap\mathcal{O}_{1}\neq\emptyset.$

(c) If $h\in\overline{W^{ss}\left(
g\right)  }$, then $\overline{W^{ss}\left(  h\right)  }\subset\overline
{W^{ss}\left(  g\right)  }.$
\end{lemma}

\begin{proof} (a) If $h\in\overline{W^{ss}\left(  c\cdot
 g\right)  }$, there exist a sequence $\left\{  h_{n}\right\}  _{n\in
\mathbb{N}}\subset W^{ss}\left(  c\cdot g\right)  $ with $h_{n}\rightarrow
h.$ It is clear from the definitions that $\left(  -c\right)  \cdot
 h_{n}\rightarrow\left(  -c\right)  \cdot h$ and $\left\{  \left(
-c\right)  \cdot h_{n}\right\}  _{n\in\mathbb{N}}\subset W^{ss}\left(
g\right)  .$ This shows that $\left(  -c\right)  \cdot h\in\overline
{W^{ss}\left(  g\right)  }$ and, hence, $h=c\cdot \bigl(\left(  -c\right)
\cdot h\bigr)\in c \cdot \left(  \overline{W^{ss}\left(  g\right)
}\right)  .$ Similarly, we show the converse inclusion.

(b) The trivialization $H:GX\rightarrow\partial^{2}X\times\mathbb{R}
$ described in section \ref{prelim} above $\bigl($see equation (\ref{lization}%
)$\bigr
)$ maps a geodesic $f\in GX$ to a triple where the third coordinate is a real
number. We will be denoting this real number by $s_{f},$ i.e.,
\begin{equation}
H\left(  f\right)  =\bigl(f\left(  -\infty\right)  ,f\left(  +\infty\right)
,s_{f}\bigr) .\label{sf}%
\end{equation}
Lift $g_{1},h_{1}$ to geodesics $\overline{g_{1}},\overline{h_{1}}\in GX$ and
consider an open neighborhood $\overline{\mathcal{O}_{1}}$ of $\overline
{h_{1}}$ of the form
\[
\overline{\mathcal{O}_{1}}:=H^{-1}\bigl(O_{1}\times O_{1}^{\prime}%
\times\left(  s_{\overline{h}_{1}}-\varepsilon_{1}^{\prime},s_{\overline
{h}_{1}}+\varepsilon_{1}^{\prime}\right)  \bigr)%
\]
where $O_{1},O_{1}^{\prime}$ are open neighborhoods of $\overline{h_{1}%
}\left(  +\infty\right)  ,$ $\overline{h_{1}}\left(  -\infty\right)  $
(respectively) in $\partial X$ with $O_{1}\cap O_{1}^{\prime}=\emptyset$ and
$\varepsilon_{1}^{\prime}$ positive real, all chosen so that
\[
p\left(  \overline{\mathcal{O}_{1}}\right)  \subseteq\mathcal{O}_{1}.
\]

\begin{claim}We may choose $\varepsilon_{1}>0$ and distinct open
neighborhoods $A_{1},A_{1}^{\prime}\subset\partial X$ of $\overline{g_{1}%
}\left(  -\infty\right)  ,\overline{g_{1}}\left(  +\infty\right)  $
(respectively) such that the neighborhood
\[
\overline{\mathcal{A}_{1}}:=H^{-1}\bigl(A_{1}\times A_{1}^{\prime}%
\times\left(  s_{\overline{g}_{1}}-\varepsilon_{1},s_{\overline{g}_{1}%
}+\varepsilon_{1}\right)  \bigr)%
\]
satisfies the following
\[
\forall\,\,\overline{g}\in\overline{\mathcal{A}_{1}}\,\,\exists\,\,\overline
{h}\in\overline{\mathcal{O}_{1}}\,\,\mathrm{such}\,\,\,\mathrm{that}%
\,\,\,\overline{h}\in W^{ss}\left(  \overline{g}\right).
\]
Then, by taking $\mathcal{A}_{1}:=p\left(  \overline{\mathcal{A}_{1}}\right)
$ the proof of the lemma is complete: for, if $g\in\mathcal{A}_{1},$ there
exists $\overline{g}\in\overline{\mathcal{A}_{1}}$ with $p\left(  \overline
{g}\right)  =g$ and, by the claim, there exists $\overline{h}\in
\overline{\mathcal{O}_{1}}$ such that $\overline{h}\in W^{ss}\left(
\overline{g}\right)  $. As $p\left(  \overline{\mathcal{O}_{1}}\right)
\subseteq\mathcal{O}_{1},$ the geodesic $h=p\left(  \overline{h}\right)  $
belongs to $\mathcal{O}_{1}$ and satisfies $h\in W^{ss}\left(  g\right) 
$. 
\end{claim}

\begin{proof}[Proof of Claim] 
Choose closed balls $B\bigl(\overline
{h_{1}}\left(  0\right)  \bigr),$ $B\bigl(\overline{g_{1}}\left(  0\right)
\bigr)
$ around $\overline{h_{1}}\left(  0\right)  ,$ $\overline{g_{1}}\left(
0\right)  $ respectively, both with radius $\varepsilon_{1}^{\prime}.$ As $X $
is proper, closed balls are compact sets and so is $\partial X.$ Thus the
(continuous) generalized Busemann function $\alpha$ restricted to $\partial
X\times B\bigl(\overline{h_{1}}\left(  0\right)  \bigr)\times B\bigl(%
\overline{g_{1}}\left(  0\right)  \bigr)$ is uniformly continuous. This
implies that for the number $\varepsilon_{1}^{\prime}/2>0,$ there exists a
compact subset $\partial B$ of $\partial X$ containing $\overline{h_{1}%
}\left(  +\infty\right)  $ and a number $\lambda>0$ such that for all $\left(
\xi,x,y\right)  \left(  \xi^{\prime},x^{\prime},y^{\prime}\right)  \in\partial
B\times B\bigl(\overline{h_{1}}\left(  0\right)  \bigr)\times B\bigl
(\overline{g_{1}}\left(  0\right)  \bigr)$ satisfying $d\left(  x,x^{\prime
}\right)  <\lambda$ and $d\left(  y,y^{\prime}\right)  <\lambda$ the following
inequality holds:
\begin{equation}
\left|  \alpha\left(  \xi,x,y\right)  -\alpha\left(  \xi^{\prime},x^{\prime
},y^{\prime}\right)  \right|  <\varepsilon_{1}^{\prime}/2 .\label{unicon}%
\end{equation}
Fix $\varepsilon_{1}<\mathrm{min}\left\{  \lambda/2,\varepsilon_{1}^{\prime
}/2\right\}  .$ We may choose small enough neighborhoods $O_{2},O_{2}^{\prime
}$ containing $\overline{h_{1}}\left(  -\infty\right)  ,\overline{h_{1}%
}\left(  +\infty\right)  $ respectively, so that if $\overline{h}$ is a
geodesic with $\overline{h}\left(  +\infty\right)  \in O_{2}^{\prime}$ and
$\overline{h}\left(  -\infty\right)  \in O_{2}$, then a suitable
re-parametrization of $\overline{h}$ (called again $\overline{h}$) satisfies
\[
d\bigl(\overline{h}\left(  0\right)  ,\overline{h_{1}}\left(  0\right)  \bigr
)<\varepsilon_{1}.
\]
We may assume that these neighborhoods $O_{2},O_{2}^{\prime}$ satisfy the
inclusions $O_{2}\subset O_{1}$ and $O_{2}^{\prime}\subset O_{1}^{\prime}%
\cap\partial B.$ Set
\[
\overline{\mathcal{O}_{2}}=H^{-1}\bigl(O_{2}\times O_{2}^{\prime}\times\left(
s_{\overline{h}_{1}}-\varepsilon_{1},s_{\overline{h}_{1}}+\varepsilon
_{1}\right)  \bigr).
\]
Then we have
\begin{equation}
\forall\,\,\overline{h}\in\overline{\mathcal{O}_{2}}\Longrightarrow d\bigl
(\overline{h}\left(  0\right)  ,\overline{h_{1}}\left(  0\right)  \bigr
)<\varepsilon_{1}+\varepsilon_{1}=2\varepsilon_{1} .\label{ho}%
\end{equation}
Moreover, using equation (\ref{flow}) and the fact that $\varepsilon
_{1}<\varepsilon_{1}^{\prime}/2$ we have
\begin{equation}
\forall\,\,s\in\left(  -\varepsilon_{1}^{\prime}/2,\varepsilon_{1}^{\prime
}/2\right)  \,\,\,and\,\,\,\forall\,\,\,\overline{h}\in\overline
{\mathcal{O}_{2}}\Longrightarrow s\cdot \overline{h}\in\overline
{\mathcal{O}_{1}} .\label{hoo}%
\end{equation}
In a similar fashion, we may choose neighborhoods $A_{1}\subset\partial X$
containing $\overline{g_{1}}\left(  -\infty\right)  $ and $A_{1}^{\prime
}\subset O_{2}^{\prime}$ containing $\overline{g_{1}}\left(  +\infty\right)
=\overline{h_{1}}\left(  +\infty\right)  $ such that for the neighborhood
\[
\overline{\mathcal{A}_{1}}=H^{-1}\bigl(A_{1}\times A_{1}^{\prime}\times\left(
s_{\overline{g}_{1}}-\varepsilon_{1},s_{\overline{g}_{1}}+\varepsilon
_{1}\right)  \bigr)%
\]
we have
\begin{equation}
\forall\,\,\,\overline{g}\in\overline{\mathcal{A}_{1}}\Longrightarrow d\bigl
(\overline{g}\left(  0\right)  ,\overline{g_{1}}\left(  0\right)  \bigr
)<2\varepsilon_{1} .\label{gio}%
\end{equation}
Let $\overline{g}\in\overline{\mathcal{A}_{1}}$ be arbitrary. Since
$A_{1}^{\prime}\subset O_{2}^{\prime},$ choose $\overline{h}\in\overline
{\mathcal{O}_{2}}$ with $\overline{h}\left(  +\infty\right)  =\overline
{g}\left(  +\infty\right)  .$ Then
\[%
\begin{array}
[c]{ll}%
d\bigl(\overline{h}\left(  0\right)  ,\overline{h_{1}}\left(  0\right)  \bigr
)<\lambda & \qquad\mathrm{by\,\,(\ref{ho})\,\,and\,\,the\,\,fact\,\,that\,}%
\text{\thinspace\thinspace}\varepsilon_{1}<\lambda/2,\\
d\bigl(\overline{g}\left(  0\right)  ,\overline{g_{1}}\left(  0\right)  \bigr
)<\lambda & \qquad\mathrm{by\,\,(\ref{gio})\,\,and\,\,the\,\,fact\,\,that\,}%
\text{\thinspace\thinspace}\varepsilon_{1}<\lambda/2,\\
\overline{h}\left(  +\infty\right)  ,\overline{h_{1}}\left(  +\infty\right)
\in\partial B & \qquad\mathrm{by\,\,construction}.
\end{array}
\]
The above three equations combined with (\ref{unicon}) imply that
\[
\left|  \alpha\bigl(\overline{h}\left(  +\infty\right)  ,\overline{h}\left(
0\right)  ,\overline{g}\left(  0\right)  \bigr)-\alpha\bigl(\overline{h_{1}%
}\left(  +\infty\right)  ,\overline{h_{1}}\left(  0\right)  ,\overline{g_{1}%
}\left(  0\right)  \bigr)\right|  <\varepsilon_{1}^{\prime}/2.
\]
As $\alpha\bigl(\overline{h}\left(  +\infty\right)  ,\overline{h_{1}}\left(
0\right)  ,\overline{g_{1}}\left(  0\right)  \bigr)=0,$ we have that
\[
-\varepsilon_{1}^{\prime}/2<\alpha\bigl(\overline{h}\left(  +\infty\right)
,\overline{h}\left(  0\right)  ,\overline{g}\left(  0\right)  \bigr
)<\varepsilon_{1}^{\prime}/2.
\]
By lemma \ref{buse}, there exists a real $s\in\left(  -\varepsilon_{1}%
^{\prime}/2,\varepsilon_{1}^{\prime}/2\right)  $ such that
\[
\alpha\bigl(\overline{g}\left(  +\infty\right)  ,s \cdot \overline{h}\left(
0\right)  ,\overline{g}\left(  0\right)  \bigr)=0.
\]
Moreover, by equation (\ref{hoo}), $s\cdot \overline{h}\in\overline
{\mathcal{O}_{1}}.$ Therefore, $s\cdot \overline{h}\in W^{ss}\left(
\overline{g}\right)  $ which completes the proof of the claim.
\renewcommand{\qed}{} 
\end{proof} 

Part (c) follows immediately from part (b).\end{proof}

\section{Topological transitivity}\label{toptransitivity}

The geodesic flow $\mathbb{R}\times GY\rightarrow GY$ is said to be
topologically transitive if given any open sets $\mathcal{O}$ and
$\mathcal{U}$ in $GY$, there exists a sequence $\left\{  t_{n}\right\}
\subset\mathbb{R},$ $t_{n}\rightarrow+\infty$ such that $t_{n}\cdot
 \mathcal{O}\,\cap\,\mathcal{U}\neq\emptyset$ for all $n\in\mathbb{N}.$ It is
apparent that topological mixing implies topological transitivity. However, in
the proof of topological mixing in section \ref{mainproof} below we will need
a property equivalent to topological transitivity, namely, that $\overline
{W^{s}\left(  f\right)  }=GY$ for any $f\in GY.$ In this section we will
establish this property without dealing with its equivalence to topological
transitivity. 

\begin{proposition}
\label{laaast}Let $X$ be a $CAT\left(  -1\right)  $-space, $\Gamma$ a
non-elementary discrete group of isometries of $X$ and $Y$ the quotient space
$Y=X$/ \negthinspace$\Gamma.$ Assume that the non-wandering set $\Omega$ of
the geodesic flow on $Y$ equals $GY.$ Then for any $f\in GY,$ $\overline
{W^{s}\left(  f\right)  }=GY.$
\end{proposition}

For the proof of the above proposition we will need the following result:

\begin{proposition}
\label{omega}Let $X,\Gamma$ and $Y$ be as above and $\Lambda\left(
\Gamma\right)  $ the limit set of the action of $\Gamma$ on $X.$ Then,
$\Omega=GY$ if and only if $\Lambda\left(  \Gamma\right)  =\partial X.$
\end{proposition}

The proof of this proposition is given in \cite{C-T} for ideal polyhedra but
it applies verbatim to our context.

\begin{proof}[Proof of Proposition \ref{laaast}] We first claim that
\begin{equation}
\mathrm{for\,\,any}\,\,\xi\in\partial X,\,\,\overline{\Gamma\xi}=\partial X.
\label{thitatrans}%
\end{equation}
Let $Fix_{h}$ be the set of points in $\partial X$ fixed by hyperbolic
elements of $\Gamma,$ i.e. $Fix_{h}=\left\{  \phi\left(  +\infty\right)
\bigm\vert\phi\in\Gamma,\phi\,\,hyperbolic\right\}  .$ As $Fix_{h}$ is dense
in $\Lambda\left(  \Gamma\right)  $ (see \cite[Ch.\ II, \S4]{Coo1}) and
$\Lambda\left(  \Gamma\right)  =\partial X$ (cf. proposition \ref{omega}) it
suffices to show that $Fix_{h}\subseteq\overline{\Gamma\xi}.$ Let $\eta\in
Fix_{h}$ be arbitrary. If $\eta=\xi$, we have nothing to show. If
$\eta\neq\xi$,
then $\eta=\phi\left(  +\infty\right)  $ for some hyperbolic $\phi\in\Gamma.$
By equation (\ref{7.2ofcoo}) it follows that $\phi^{n}\left(  \xi\right)
\rightarrow\eta,$ hence, $\eta\in\overline{\Gamma\xi}.$ This completes the
proof of equation (\ref{thitatrans}).

Now let $f,g\in GY$ be
arbitrary. We proceed to find a sequence $f_{n}\in W^{s}\left(  f\right)  $
such that $f_{n}\rightarrow g.$ Lift $f,g$ to geodesics $\overline
{f},\overline{g}$ in $GX$ and let $s_{\overline{g}}$ be the unique real number
so that $\bigl(\overline{g}\left(  -\infty\right)  ,\overline{g}\left(
+\infty\right)  ,s_{\overline{g}}\bigr)=H\left(  \overline{g}\right)  $ (cf.
equation \ref{sf}). By equation (\ref{thitatrans}) there exists a sequence
$\left\{  \phi_{n}\right\}  \subset\Gamma$ such that $\phi_{n}\bigl
(\overline{f}\left(  +\infty\right)  \bigr)\rightarrow\overline{g}\left(
+\infty\right)  .$ Define a sequence $\left\{  \overline{g_{n}}\right\}
\subset GX$ where each $\overline{g_{n}}$ is determined by the following three
conditions (cf. equation \ref{lization}):
\[%
\begin{array}
[c]{l}%
\mathrm{(i)\,\,\,}\overline{g_{n}}\left(  +\infty\right)  =\phi_{n}\bigl
(\overline{f}\left(  +\infty\right)  \bigr),\\
\mathrm{(ii)\,\,\,}\overline{g_{n}}\left(  -\infty\right)
\,\,\mathrm{is\,\,any\,\,sequence}\,:\overline{g_{n}}\left(  -\infty\right)
\rightarrow\overline{g}\left(  -\infty\right)  ,\\
\mathrm{(iii)\,\,\,the\,\,parametrization\,\,is\,\,chosen\,\,by\,\,requiring\,}%
\,s_{\overline{g_{n}}}=s_{\overline{g}}\,\,\mathrm{\forall}\,\,n.
\end{array}
\]
In other words, $\overline{g_{n}}:=H^{-1}\Bigl(\overline{g}\left(
-\infty\right)  ,\phi_{n}\bigl(\overline{f}\left(  +\infty\right)  \bigr
),s_{\overline{g}}\Bigr).$ It is apparent that $\overline{g_{n}}%
\rightarrow\overline{g}.$ Define $\overline{f_{n}}:=\phi_{n}^{-1}\left(
\overline{g_{n}}\right)  $ and set $f_{n}=p\left(  \overline{f_{n}}\right)  .$
As $p\left(  \overline{f_{n}}\right)  =p\left(  \overline{g_{n}}\right)  $ and
$\overline{g_{n}}\rightarrow\overline{g},$ it follows that $f_{n}\rightarrow
p\left(  \overline{g}\right)  =g.$ Moreover,
\[
\overline{f_{n}}\left(  +\infty\right)  =\phi_{n}^{-1}\bigl(\overline
{g}\left(  +\infty\right)  \bigr)=\phi_{n}^{-1}\Bigl(\phi_{n}\bigl(%
\overline{f}\left(  +\infty\right)  \bigr)\Bigr)=\overline{f}\left(
+\infty\right)
\]
thus, $f_{n}\in W^{s}\left(  f\right)  $ as required.\end{proof}

\begin{corollary}
\label{foranyclosed}Under the assumptions of proposition \ref{laaast} above,
if there exists a geodesic $g$ whose strong stable set $W^{ss}\left(
g\right)  $ satisfies $\overline{W^{ss}\left(  g\right)  }=GY$, then
$\overline{W^{ss}\left(  f\right)  }=GY$ for any closed geodesic $f\in GY.$
\end{corollary}

\begin{proof}Let $g$ be a geodesic in $GY$ satisfying
$\overline{W^{ss}\left(  g\right)  }=GY$ and let $f$ be a closed geodesic in
$GY$ with period, say, $\omega.$ By proposition \ref{laaast}, $\overline
{W^{s}\left(  f\right)  }=GY$ hence, $g\in\overline{W^{s}\left(  f\right)  }.$
This means that there exists a sequence $\left\{  g_{n}\right\}  \subset
W^{s}\left(  f\right)  $ such that $g_{n}\rightarrow g.$ For each
$n\in\mathbb{N},$ consider lifts $\overline{g_{n}},\overline{f}$ of $g_{n},f$
respectively, satisfying $\overline{g_{n}}\in W^{ss}\left(  \overline
{f}\right)  $ and use lemma \ref{buse}(b) to obtain a real number $t_{n}$ such
that $t_{n}\cdot g_{n}\in W^{ss}\left(  f\right)  .$ Each $t_{n}$ may be
expressed by 
\[
t_{n}=k\omega+c_{n}%
\]
where $k\in\mathbb{Z}$ and $c_{n}\in\left[  0,\omega\right)  .$ By choosing,
if necessary, a subsequence we have that $c_{n}\rightarrow c$ for some
$c\in\left[  0,\omega\right]  .$ Then $c_{n}\cdot g_{n}\rightarrow
c\cdot g$ with $c_{n}\cdot g_{n}\in W^{ss}\left(  f\right)  .$ This
means that for some $c\in\left[  0,\omega\right]  ,$ $c \cdot g\in
\overline{W^{ss}\left(  f\right)  }.$ By lemma \ref{contofsets}(a), we have
$c\cdot \overline{W^{ss}\left(  g\right)  }=\overline{W^{ss}\left(
c\cdot g\right)  }$ and, by lemma \ref{contofsets}(c), $\overline
{W^{ss}\left(  c\cdot g\right)  }\subset\overline{W^{ss}\left(  f\right)
}.$ Thus, $GY=c\cdot GY=c\cdot \overline{W^{ss}\left(  g\right)
}\subset\overline{W^{ss}\left(  f\right)  }.$\end{proof}

\section{Proof of topological mixing}\label{mainproof}

For the proof of theorem \ref{main} we follow closely the idea used by
Eberlein in \cite{Ebe2} in the proof of topological mixing of the geodesic
flow on Riemannian manifolds of non-positive curvature. However, since we deal
with a more general class of spaces, the difficulties which arise here are of
a different nature. We first establish the following:

\begin{proposition}
\label{ss}Let $X,Y$ and $\Gamma$ be as in theorem \ref{main} above. If
$\Omega=GY$, then there exists a geodesic $g$ whose strong stable set
$W^{ss}\left(  g\right)  $ satisfies $\overline{W^{ss}\left(  g\right)  }=GY.$
\end{proposition}

\begin{proof} We first show that
\begin{equation}
\forall\,\,\mathcal{O},\mathcal{U}\subseteq GY\,\,open,\,\exists
\,\,g\in\mathcal{O}:W^{ss}\left(  g\right)  \cap\mathcal{U}\neq\emptyset.
\label{arguem}%
\end{equation}
We will use the letter $p$ to denote both projections $X\rightarrow Y$ and
$GX\rightarrow GY.$ Let $\mathcal{O},\mathcal{U}\subseteq GY\,$be arbitrary
open sets. Using equation (\ref{5.1ofcoo}) and the fact that $\partial
X=\Lambda\left(  \Gamma\right)  $ (which follows from the assumption
$\Omega=GY $ and proposition \ref{omega} above), we may choose $f\in
p^{-1}\left(  \mathcal{O}\right)  $ and $h\in p^{-1}\left(  \mathcal{U}%
\right)  $ such that
\[
\bigl(f\left(  +\infty\right)  ,h\left(  +\infty\right)  \bigr)=\left(
\phi\left(  +\infty\right)  ,\phi\left(  -\infty\right)  \right)  :\phi
\in\Gamma\,\,\,is\,\,\,hyperbolic.
\]
For each $n,$ let $\xi_{n}$ be in $\partial X$ such that $\alpha\bigl(\xi
_{n},f\left(  0\right)  ,\phi^{n}\left(  h\left(  0\right)  \right)  \bigr
)=0.$ We claim that
\begin{equation}
\xi_{n}\rightarrow f\left(  +\infty\right)  \,\,as\,\,n\rightarrow\infty.
\label{xn}%
\end{equation}
To see this assume, on the contrary, that $\left\{  \xi_{n}\right\}  $ (or, a
subsequence of it) converges to $\xi\in\partial X$ with $\xi\neq f\left(
+\infty\right)  .$ Let $\beta$ be a geodesic in $X$ such that $\beta\left(
+\infty\right)  =f\left(  +\infty\right)  $ and $\beta\left(  -\infty\right)
=\xi.$ Similarly, let $\beta^{\prime}\in GX$ such that $\beta^{\prime}\left(
+\infty\right)  =f\left(  +\infty\right)  $ and $\beta^{\prime}\left(
-\infty\right)  =h\left(  +\infty\right)  .$ Since $\phi$ translates
$\beta^{\prime},$ i.e., $\phi\bigl(\mathrm{Im}\left(  \beta^{\prime}\right)
\bigr)=\mathrm{Im}\left(  \beta^{\prime}\right)  ,$ it follows that for any
$n\in\mathbb{N}$%
\[
\mathrm{dist}\bigl(\phi^{n}\left(  h\left(  0\right)  \right)  ,\mathrm{Im}%
\beta^{\prime}\bigr)\leq\mathrm{dist}\bigl(h\left(  0\right)  ,\mathrm{Im}%
\beta^{\prime}\bigr).
\]
As $\beta$ and $\beta^{\prime}$ are asymptotic and $\phi^{n}\left(  h\left(
0\right)  \right)  \rightarrow\beta\left(  +\infty\right)  =\beta^{\prime
}\left(  +\infty\right)  $, a similar statement holds true for $\beta$ by
using condition (C), namely,
\[
\exists\,M\in\mathbb{R}:\mathrm{dist}\bigl(\phi^{n}\left(  h\left(  0\right)
\right)  ,\mathrm{Im}\beta\bigr)\leq M,\,\,for\,\,all\,\,n\in\mathbb{N}.
\]
For each $n\in\mathbb{N},$ let $t_{n}$ be the real number realizing the
distance in the left-hand side of the above equation (cf. lemma \ref{pro}). By
equation (\ref{7.2ofcoo}), is clear that $t_{n}\rightarrow+\infty.$ As the
generalized Busemann function $\alpha$ is Lipschitz (with Lipschitz constant
1) with respect to the third variable, it follows that for all $n\in
\mathbb{N},$%
\begin{equation}
\left|  \alpha\bigl(\xi,f\left(  0\right)  ,\beta\left(  t_{n}\right)  \bigr
)-\alpha\bigl(\xi,f\left(  0\right)  ,\phi^{n}\left(  h\left(  0\right)
\right)  \bigr)\right|  \leq M .\label{zero}%
\end{equation}
$\xi_{n}$ is chosen so that $\phi^{n}\left(  h\left(  0\right)  \right)  $ and
$f\left(  0\right)  $ are equidistant from $\xi_{n},$ thus,
\[
\mathrm{lim}_{n\rightarrow\infty}\alpha\bigl(\xi,f\left(  0\right)  ,\phi
^{n}\left(  h\left(  0\right)  \right)  \bigr)=\mathrm{lim}_{n\rightarrow
\infty}\alpha\bigl(\xi_{n},f\left(  0\right)  ,\phi^{n}\left(  h\left(
0\right)  \right)  \bigr)=0.
\]
The latter combined with equation (\ref{zero}) implies that
\[
\left|  \mathrm{lim}_{n\rightarrow\infty}\alpha\bigl(\xi,f\left(  0\right)
,\beta\left(  t_{n}\right)  \bigr)\right|  \leq M.
\]
This is impossible by lemma \ref{buse}(c) and the fact that $t_{n}%
\rightarrow+\infty.$ Thus equation (\ref{xn}) is proved. We next show that
\begin{equation}
\phi^{-n}\left(  \xi_{n}\right)  \rightarrow h\left(  +\infty\right)
\,\,as\,\,n\rightarrow\infty. \label{yn}%
\end{equation}
Assume, on the contrary, that $\phi^{-n}\left(  \xi_{n}\right)  $ (or, a
subsequence of it) converges to $\zeta\in\partial X$ with $\zeta\neq h\left(
+\infty\right)  .$ We choose a geodesic $\beta\in GX$ with $\beta\left(
+\infty\right)  =\zeta$ and $\beta\left(  -\infty\right)  =h\left(
+\infty\right)  $ and proceed with the proof exactly as in the previous
argument by using the facts that
\[
\mathrm{lim}_{n\rightarrow\infty}\phi^{-n}\left(  f\left(  0\right)  \right)
=h\left(  +\infty\right)  
\]
and
\[
\alpha\bigl(\phi^{-n}\left(  \xi_{n}\right)  ,\phi^{-n}\left(  f\left(
0\right)  \right)  ,h\left(  0\right)  \bigr)=0.
\]
Choose now geodesics $f_{n}\in GX,$ $n\in\mathbb{N}$ such that $f_{n}\left(
+\infty\right)  =\xi_{n}$ and $f_{n}\left(  -\infty\right)  =f\left(
-\infty\right)  .$ We may parametrize $f_{n}$ so that $f_{n}\left(  0\right)
\rightarrow f\left(  0\right)  .$ This can be done by requiring $s_{\overline
{f_{n}}}=s_{\overline{f}}$ for all $n\in\mathbb{N}.$ $\bigl($cf. equation
(\ref{sf})$\bigr)$. Similarly, choose $h_{n}\in GX$ such that $h_{n}\left(
+\infty\right)  =\xi_{n}$ and $h_{n}\left(  -\infty\right)  =\phi^{n}\bigl
(h\left(  -\infty\right)  \bigr)$ and parametrize them so that
\begin{equation}
\alpha\bigl(\xi_{n},f_{n}\left(  0\right)  ,h_{n}\left(  0\right)  \bigr)=0.
\label{finalb}%
\end{equation}
It is apparent that for $n$ large enough, $f_{n}\in p^{-1}\left(
\mathcal{O}\right)  $ and $h_{n}\in W^{ss}\left(  f_{n}\right)  $. If we show
that $\phi^{-n}\left(  h_{n}\right)  \in p^{-1}\left(  \mathcal{U}\right)  $
for $n$ large enough, then we would have
\[%
\begin{array}
[c]{l}%
p\left(  f_{n}\right)  \in\mathcal{O},\\
p\left(  h_{n}\right)  =p\bigl(\phi^{-n}\left(  h_{n}\right)  \bigr)%
\in\mathcal{U},\\
p\left(  h_{n}\right)  \in W^{ss}\bigl(p\left(  f_{n}\right)  \bigr).
\end{array}
\]
The above three properties imply that for $n$ large enough, $W^{ss}\bigl
(p\left(  f_{n}\right)  \bigr)\cap\mathcal{U}\neq\emptyset,$ as required in
equation (\ref{arguem}). We conclude the proof of the proposition by showing
that $\phi^{-n}\left(  h_{n}\right)  \in p^{-1}\left(  \mathcal{U}\right)  .$
Using equation (\ref{yn}) above, it is clear that
\begin{equation}%
\begin{array}
[c]{lll}%
\bigl(\phi^{-n}\left(  h_{n}\right)  \bigr)\left(  +\infty\right)  & = &
\phi^{-n}\bigl(h_{n}\left(  \infty\right)  \bigr)\\
& = & \phi^{-n}\left(  \xi_{n}\right)  \rightarrow h\left(  +\infty\right).
\end{array}
\label{a}%
\end{equation}
Similarly,
\begin{equation}%
\begin{array}
[c]{lll}%
\bigl(\phi^{-n}\left(  h_{n}\right)  \bigr)\left(  -\infty\right)  & = &
\phi^{-n}\bigl(h_{n}\left(  -\infty\right)  \bigr)\\
& = & \phi^{-n}\bigl(\phi^{n}\bigl(h\left(  -\infty\right)  \bigr)\bigr
)=h\left(  -\infty\right)  \,\,as\,\,n\rightarrow\infty.
\end{array}
\label{c}%
\end{equation}
By condition (U) and equations (\ref{a}), (\ref{c}), we have
\[
d\Bigl(h\left(  0\right)  ,\mathrm{Im}\,\phi^{-n}\left(  h_{n}\right)
\Bigr)\rightarrow0
\]
as $n\rightarrow+\infty,$ and, therefore,
\begin{equation}
d\Bigl(\phi^{n}\bigl(h\left(  0\right)  \bigr),\mathrm{Im}\,h_{n}%
\Bigr)\rightarrow0\,\,as\,\,n\rightarrow\infty. \label{finala}%
\end{equation}
Let $h_{n}\left(  t_{n}\right)  ,t_{n}\in\mathbb{R}$ be the point on
$\mathrm{Im}\,h_{n}$ which realizes the distance in equation (\ref{finala})
above (cf. lemma \ref{pro}). As the function $\alpha$ is Lipschitz with
respect to the third variable (with Lipschitz constant 1), we have
\[
\left|  \alpha\bigl(\xi_{n},f\left(  0\right)  ,\phi^{n}\bigl(h\left(
0\right)  \bigr)\bigr)-\alpha\bigl(\xi_{n},f\left(  0\right)  ,h_{n}\left(
t_{n}\right)  \bigr)\right|  \leq d\Bigl(\phi^{n}\bigl(h\left(  0\right)
\bigr),h_{n}\left(  t_{n}\right)  \Bigr).
\]
Using the defining property of $\xi_{n},$ i.e., $\alpha\bigl(\xi_{n},f\left(
0\right)  ,\phi^{n}\left(  h\left(  0\right)  \right)  \bigr)=0,$ it follows
that
\[
\alpha\bigl(\xi_{n},f\left(  0\right)  ,h_{n}\left(  t_{n}\right)  \bigr
)\rightarrow0\,\,as\,\,n\rightarrow\infty.
\]
Similarly, using the fact that $f_{n}\left(  0\right)  \rightarrow f\left(
0\right)  $ as $n\rightarrow\infty$ and the Lipschitz property of $\alpha$
with respect to the second variable, we have
\[
\alpha\bigl(\xi_{n},f_{n}\left(  0\right)  ,h_{n}\left(  t_{n}\right)  \bigr
)\rightarrow0\,\,as\,\,n\rightarrow\infty.
\]
Since, by lemma \ref{buse}(b), there is a unique point on each $\mathrm{Im}%
\,h_{n}$ which is equidistant from $f_{n}\left(  0\right)  $ with respect to
$\xi_{n},$ namely, $h_{n}\left(  0\right)  $ $\bigl($cf. equation
(\ref{finalb})$\bigr),$ it follows that $t_{n}\rightarrow0$ which, combined
with equation (\ref{finala}) implies that
\[
d\Bigl(\phi^{n}\bigl(h\left(  0\right)  \bigr),h_{n}\left(  0\right)
\Bigr)\rightarrow0\,\,as\,\,n\rightarrow\infty.
\]
Therefore, $\phi^{-n}\bigl(h_{n}\left(  0\right)  \bigr)\rightarrow h\left(
0\right)  $ as $n\rightarrow\infty.$ The latter combined with the facts
\begin{align*}
\bigl(\phi^{-n}\left(  h_{n}\right)  \bigr)\left(  +\infty\right)   &
\rightarrow h\left(  +\infty\right)  \,\,as\,\,n\rightarrow\infty,\\
\bigl(\phi^{-n}\left(  h_{n}\right)  \bigr)\left(  -\infty\right)   &
\rightarrow h\left(  -\infty\right)  \,\,as\,\,n\rightarrow\infty
\end{align*}
implies that $\phi^{-n}\left(  h_{n}\right)  \in p^{-1}\left(  \mathcal{U}%
\right)  $ which concludes the proof of equation (\ref{arguem}).

Using now a countable basis $\left\{  \mathcal{O}_{n}\right\}  _{n\in
\mathbb{N}}$ for the topology of $GY,$ which exists by equation
(\ref{secondcount}), the proof is completed by a standard topological argument
(cf. \cite[Th.\ 5.2]{Ebe2}) which we include here for the readers convenience:
If $\mathcal{O}$ is an arbitrary open set in $GY$, then by equation
(\ref{arguem}) above, there exists $g_{1}\in\mathcal{O}$ such that
$W^{ss}\left(  g_{1}\right)  \cap\mathcal{O}_{1}\neq\emptyset.$ Let $h_{1}\in
W^{ss}\left(  g_{1}\right)  \cap\mathcal{O}_{1}.$ By lemma \ref{contofsets},
there exists an open set $\mathcal{A}_{1}$ containing $g_{1}$ satisfying
\[
W^{ss}\left(  h\right)  \cap\mathcal{O}_{1}\neq\emptyset
\]
for every $h\in\mathcal{A}_{1}.$ Moreover, we may assume that the closure
$\overline{\mathcal{A}_{1}}$ of $\mathcal{A}_{1}$ lies in $\mathcal{O}_{1}$
and is compact. Inductively, a sequence of open sets $\mathcal{A}_{i}$ and a
sequence of geodesics $g_{i}$ are constructed such that
\begin{align*}
&\overline{\mathcal{A}_{i}}\subset\mathcal{A}_{i-1},\\
&g_{i}\in\mathcal{A}_{i},\\
&for\,\,every\,\,h\in\mathcal{A}_{i},W^{ss}\left(  h\right)  \cap
\mathcal{O}_{i}\neq\emptyset.
\end{align*}
By the finite intersection property of the compact sets $\overline
{\mathcal{A}_{i}}$ it follows that there exists a $g\in\bigcap_{i=1}^{\infty
}\overline{\mathcal{A}_{i}}.$ As $W^{ss}\left(  g\right)  \cap\mathcal{O}%
_{i}\neq\emptyset$ for all $i$, $\overline{W^{ss}\left(  g\right)  }%
=GY.$\end{proof}

We will need a pointwise version of topological mixing and a criterion for
such property.

\begin{definition}
Let $h,f$ be in $GY$ and let $\left\{  s_{n}\right\}  _{n\in\mathbb{N}}$ be a
sequence converging to $+\infty$ or $-\infty.$ We say that $h$ is $s_{n}%
$-mixing with $f$ (notation, $h\sim_{s_{n}}f$) if for every neighborhood 
$\mathcal{O},\mathcal{U}$ in $GY$ of $h,f$ respectively, $s_{n}\cdot
 \mathcal{O}\cap\mathcal{U}\neq\emptyset$ for all $n$ sufficiently large$.$
\end{definition}

If $h\sim_{s_{n}}f$ for some $h,f\in GY$, then using decreasing
sequences of open neighborhoods of $h$ and $f$ it is easily shown that for
each subsequence $\left\{  s_{n}^{\prime}\right\}  $ of $\left\{
s_{n}\right\}  $ there exists a subsequence $\left\{  r_{n}\right\}  $ of
$\left\{  s_{n}^{\prime}\right\}  $ and a sequence $\left\{  h_{n}\right\}
\subset GY$ such that $h_{n}\rightarrow h$ and $r_{n}\cdot h_{n}%
\rightarrow f.$ The proof of the converse statement is elementary, hence, the
following criterion for the $s_{n}$-mixing of $h,f$ holds.

\begin{criterion}
\label{criterion}If $h,f\in GY$, then $h\sim_{s_{n}}f$ if and only if for each
subsequence $\left\{  s_{n}^{\prime}\right\}  $ of $\left\{  s_{n}\right\}  $
there exists a subsequence $\left\{  r_{n}\right\}  $ of $\left\{
s_{n}^{\prime}\right\}  $ and a sequence $\left\{  h_{n}\right\}  \subset GY$
such that $h_{n}\rightarrow h$ and $r_{n}\cdot h_{n}\rightarrow f.$
\end{criterion}

The following lemma asserts that pointwise topological mixing is
transferred via the strong stable relation of geodesics.

\begin{lemma}
\label{new}If $f,g,g^{\prime}\in GY$ so that $f\in\overline{W^{ss}\left(
g\right)  }$ and $g\sim_{s_{n}}g^{\prime}$ for some sequence $s_{n}%
\rightarrow\infty,$ then $f\sim_{s_{n}}g^{\prime}.$
\end{lemma}

\begin{proof} Fix a sequence $\left\{  s_{n}\right\}  $ with
$s_{n}\rightarrow\infty.$ It is easy to verify that the set
\[
\left\{  f\in GY\bigm\vert f\sim_{s_{n}}g^{\prime}\right\}
\]
is closed; for, if $\left\{  h_{k}\right\}  _{k\in\mathbb{N}}\subset GY$ with
$h_{k}\rightarrow h$ and each $h_{k}$ is $s_{n}$-mixing with $g^{\prime},$ let
$\mathcal{O},\mathcal{U}$ be neighborhoods of $h,g^{\prime}$ respectively.
Since $h_{k}\rightarrow h,$ $\mathcal{O}$ is also a neighborhood of $h_{k_{0}%
},$ for some $k_{0}$ large enough. As $h_{k_{0}}\sim_{s_{n}}g^{\prime}$, it
follows that $s_{n}\cdot \mathcal{O}\cap\mathcal{U}\neq\emptyset$ for $n$
sufficiently large which implies that $\left\{  f\in GY\bigm\vert f\sim
_{s_{n}}g^{\prime}\right\}  $ is closed. Therefore, it suffices to prove the
assertion of the lemma for $f\in W^{ss}\left(  g\right)  $.

In order to use criterion \ref{criterion} above to show that $f\sim
_{s_{n}}g^{\prime},$ let $\left\{  t_{n}\right\}  $ be an arbitrary subsequence
of $\left\{  s_{n}\right\}  .$ As $g\sim_{s_{n}}g^{\prime}$, there exists
(again by criterion \ref{criterion}) a subsequence $\left\{  r_{n}\right\}  $
of $\left\{  t_{n}\right\}  $ and a sequence $\left\{  g_{n}\right\}  $
converging to $g$ such that $r_{n}\cdot g_{n}\rightarrow g^{\prime}.$
Lift $g$ and $f$ to geodesics $\overline{g}$ and $\overline{f}$ in $GX$ such
that $\overline{f}\left(  +\infty\right)  =\overline{g}\left(  +\infty\right)
$ and $\alpha\left(  \overline{f}\left(  +\infty\right)  ,\overline{f}\left(
0\right)  ,\overline{g}\left(  0\right)  \right)  =0$. Lift each $g_{n}$ to a
geodesic $\overline{g_{n}}$ such that $\overline{g_{n}}\left(  +\infty\right)
\rightarrow\overline{g}\left(  +\infty\right)  ,$ $\overline{g_{n}}\left(
-\infty\right)  \rightarrow\overline{g}\left(  -\infty\right)  $ and
$\overline{g_{n}}\left(  0\right)  \rightarrow\overline{g}\left(  0\right)  $.
Define a sequence of geodesics $\left\{  \overline{f_{n}}\right\}
_{n\in\mathbb{N}}$ such that $\overline{f_{n}}\rightarrow\overline{f}$ with
$\overline{f_{n}}\left(  +\infty\right)  =\overline{g_{n}}\left(
+\infty\right)  $ and $\overline{f_{n}}\left(  -\infty\right)  =\overline
{f}\left(  -\infty\right)  .$ By the continuity of the $\alpha$ function we
have that
\[
\mathrm{lim}_{n\rightarrow\infty}\alpha\bigl(\xi_{n},\overline{f_{n}}\left(
0\right)  ,\overline{g_{n}}\left(  0\right)  \bigr)=\alpha\bigl(\xi
,\overline{f}\left(  0\right)  ,\overline{g}\left(  0\right)  \bigr)=0
\]
hence, by passing if necessary to a subsequence of $\left\{  \overline{f_{n}%
}\right\}  _{n\in\mathbb{N}},$ we may assume that
\[
\alpha\bigl(\xi_{n},\overline{f_{n}}\left(  0\right)  ,\overline{g_{n}}\left(
0\right)  \bigr)<1/n,\,\,\,for\,\,all\,\,n\in\mathbb{N}.
\]
By lemma \ref{buse}(b) we may choose the parametrization of each
$\overline{f_{n}}$ so that
\begin{equation}
\alpha\bigl(\xi_{n},\overline{f_{n}}\left(  0\right)  ,\overline{g_{n}}\left(
0\right)  \bigr)=0,\,\,\,for\,\,all\,\,n\in\mathbb{N} .\label{efen}%
\end{equation}
As the change of parametrization tends to $0$ as $n\rightarrow\infty$, we may
assume that the sequence $\left\{  \overline{f_{n}}\right\}  _{n\in\mathbb{N}%
}$ satisfies equation (\ref{efen}) and $\overline{f_{n}}\rightarrow
\overline{f}.$ Moreover, if we set $\,f_{n}:=p\left(  \overline{f_{n}}\right)
$, then $f_{n}\rightarrow f.$

We proceed now to show that $r_{n}\cdot f_{n}\rightarrow g^{\prime}.$ Let
$K$ be an arbitrary compact subset of $\mathbb{R}$ and $\varepsilon$ arbitrary
positive. By construction, $\overline{f_{n}}\in W^{ss}\left(  \overline{g_{n}%
}\right)  $ for all $n\in\mathbb{N}$ and $\overline{f}\in W^{ss}\left(
\overline{g}\right)  .$ Moreover, by proposition \ref{asy},
\begin{equation}%
\begin{array}
[c]{l}%
\mathrm{lim}_{t\rightarrow\infty}d\bigl(\overline{f_{n}}\left(  t\right)
,\overline{g_{n}}\left(  t\right)  \bigr)=0,\\
\mathrm{lim}_{t\rightarrow\infty}d\bigl(\overline{f}\left(  t\right)
,\overline{g}\left(  t\right)  \bigr)=0.
\end{array}
\label{double}%
\end{equation}
Choose a positive real $T$ such that
\[
d\bigl(\overline{f}\left(  T\right)  ,\overline{g}\left(  T\right)  \bigr
)<\varepsilon/6.
\]
The above equation holds for all $t>T.$ This follows by convexity of the
distance function (see \cite[Ch. 2]{B-G-H}) and equation (\ref{double}). As
$\overline{f_{n}}\rightarrow\overline{f}$ and $\overline{g_{n}}\rightarrow
\overline{g}$ we may choose $N\in\mathbb{N}$ such that
\[%
\begin{array}
[c]{l}%
d\bigl(\overline{f_{n}}\left(  T\right)  ,\overline{f}\left(  T\right)  \bigr
)<\varepsilon/6,\\
d\bigl(\overline{g_{n}}\left(  T\right)  ,\overline{g}\left(  T\right)  \bigr
)<\varepsilon/6.
\end{array}
\]
Thus, $d\bigl(\overline{f_{n}}\left(  T\right)  ,\overline{g_{n}}\left(
T\right)  \bigr
)<\varepsilon/2$ and as before, it follows that
\[
d\bigl(\overline{f_{n}}\left(  t\right)  ,\overline{g_{n}}\left(  t\right)
\bigr)<\varepsilon/2\,\,\,\mathrm{for\,\,\,all\,\,}t>T.
\]
As $r_{n}\rightarrow+\infty,$ there exists $n_{0}$ such that $r_{n}\geq
T+diam\,K$ for all $n\geq n_{0}.$ Now for all $n$ sufficiently large, namely,
$n\geq\mathrm{max}\left\{  N,n_{0}\right\}  ,$ we have
\[
d\bigl(\overline{f_{n}}\left(  r_{n}+t\right)  ,\overline{g_{n}}\left(
r_{n}+t\right)  \bigr)<\varepsilon/2,\,\,\forall\,\,t\in K
\]
which implies that
\[
d\bigl(r_{n}\cdot f_{n}\left(  t\right)  ,r_{n}\cdot g_{n}\left(
t\right)  \bigr)<\varepsilon/2,\,\,\forall\,\,t\in K.
\]
As $r_{n}\cdot g_{n}\rightarrow g^{\prime},$ we have that for all $n$
sufficiently large
\[
d\bigl(r_{n}\cdot g_{n}\left(  t\right)  ,g^{\prime}\left(  t\right)
\bigr)<\varepsilon/2,\,\,\forall\,\,t\in K.
\]
Combining the last two inequalities we obtain that
\[
d\bigl(r_{n}\cdot f_{n}\left(  t\right)  ,g^{\prime}\left(  t\right)
\bigr)<\varepsilon,\,\,\forall\,\,t\in K.
\]
As $K,\varepsilon$ were arbitrary, we have shown that for all $n$ sufficiently
large, $r_{n}\cdot f_{n}$ lies in any neighborhood of $g^{\prime}.$
Therefore, $r_{n}\cdot f_{n}\rightarrow g^{\prime}$ as required.\end{proof}

\begin{proof}[Proof of Theorem \ref{main}] We will show the following
property:
\begin{equation}
\forall\,\,\mathcal{O},\mathcal{U}\subseteq GY\,\,open,\ \exists\,\,t_{0}%
^{\prime}>0:t\cdot \mathcal{O}\cap\mathcal{U}\neq\emptyset\,\,\forall
\,\,t\geq t_{0}^{\prime} .\label{suffices}%
\end{equation}
Then applying this property to the open sets $-\mathcal{O},-\mathcal{U}$ we
obtain a number
\begin{equation}
t_{0}^{\prime\prime}>0:\bigl(t\cdot\left(  -\mathcal{O}\right)  \bigr)%
\cap\left(  -\mathcal{U}\right)  \neq\emptyset\,\,\forall\,\,t\geq
t_{0}^{\prime\prime}. \label{suffices-}%
\end{equation}
Setting $t_{0}=\max\left\{  t_{0}^{^{\prime}},t_{0}^{^{\prime\prime}}\right\}
$ we have that
\[
-\bigl(\left(  -t\right)  \cdot \mathcal{O}\cap\mathcal{U}\bigr)=\bigl
(-\left(  \left(  -t\right)  \cdot \mathcal{O}\right)  \bigr)\cap\left(
-\mathcal{U}\right)  =\bigl(t\cdot \left(  -\mathcal{O}\right)  \bigr)%
\cap\left(  -\mathcal{U}\right)  .
\]
By (\ref{suffices-}) it follows that $-\bigl(\left(  -t\right)  \cdot
 \mathcal{O}\cap\mathcal{U}\bigr)\neq\emptyset,$ hence $t\cdot
 \mathcal{O}\cap\mathcal{U}\neq\emptyset$ for all $t\leq-t_{0}.$ This
combined with (\ref{suffices}) completes the proof of theorem \ref{main}. We
proceed now to show equation (\ref{suffices}). For this it suffices to show
that
\begin{equation}%
\begin{array}
[c]{l}%
\forall\,h,f\in GY\,\,\mathrm{and}\,\,\forall\,\left\{  t_{n}\right\}
\,\mathrm{with}\,t_{n}\rightarrow\infty,\exists\,\,\mathrm{sub}\text{-}\\
\mathrm{sequence}\,\,\left\{  s_{n}\right\}  \subset\left\{  t_{n}\right\}
\,\,\mathrm{such}\,\,\mathrm{that}\,\,h\sim_{s_{n}}f.
\end{array}
\label{sigma}%
\end{equation}
Let $f,h$ and $\left\{  t_{n}\right\}  $ be given and $g$ be the geodesic
provided by proposition \ref{ss}. By corollary \ref{foranyclosed} we may
assume that $g$ is closed. Moreover, by lemma \ref{contofsets}(a), the
conclusion of proposition \ref{ss} is satisfied by any translate $c\cdot g
$ of $g$, where $c\in\mathbb{R}.$ Choose a subsequence $\left\{
s_{n}\right\}  $ of $\left\{  t_{n}\right\}  $ such that $s_{n}\cdot
 g\rightarrow c\cdot g$ for some $c\in\left[  0,\textup{period}\left(  g\right)
\right]  .$ It is apparent that $g\sim_{s_{n}}c\cdot g.$ Since
$\overline{W^{ss}\left(  g\right)  }=GY,$ $f\in\overline{W^{ss}\left(
g\right)  }$ and hence, by lemma \ref{new}, $f\sim_{s_{n}}c\cdot g.$ This
implies that $-c\cdot g\sim_{s_{n}}-f.$ Applying lemma \ref{new} again and
using the fact that $-h\in\overline{W^{ss}\left(  -c\cdot g\right)  }$ it
follows that $-h\sim_{s_{n}}-f,$ thus $f\sim_{s_{n}}h$ as required.\end{proof}

\section{Applications}\label{application}

In this section we provide classes of spaces, much wider than Riemannian
manifolds, satisfying all assumptions posited in theorem \ref{main} above.
Recall that a metric space is \textit{geodesically complete} if each geodesic
segment is the restriction of a geodesic defined on the whole real line. An
immediate application is the following

\begin{corollary}
\label{Bestvina}Let $X$ be a proper geodesically complete $CAT\left(
-1\right)  $-space and $\Gamma$ a discrete one-ended group of isometries of $X$ 
with compact quotient $Y=X/\Gamma.$ Then the geodesic flow on $Y$ is
topologically mixing.
\end{corollary}

\begin{proof} $X$ is the universal cover of $Y$ and a hyperbolic
space in the sense of Gromov. Since $Y$ is compact, $\Gamma\approx\pi
_{1}\left(  Y\right)  $ is a hyperbolic group whose boundary is isomorphic
with the boundary of $X$ (see \cite[Ch.\ 4, Theorem\ 4.1]{C-D-P}). In particular,
$\Gamma$ is non-elementary. Since $\Gamma$ is one-ended, it follows that the
boundary $\partial X$ of $X$ is connected (see for example Bowditch
\cite{Bow1,Bow2}). It is easy now to deduce that for any $x,x^{\prime}\in X$
there exists a point $\xi\in\partial X$ such that $\alpha\left(
\xi,x,x^{\prime}\right)  =0$: if $x=x^{\prime}$, the result is trivial. If
$x\neq x^{\prime},$ extend the geodesic segment joining $x$ with $x^{\prime}$
to a geodesic, say, $g.$ We may assume that $x=g\left(  s\right)  ,$
$x^{\prime}=g\left(  s^{\prime}\right)  $ for some $s,s^{\prime}\in\mathbb{R}$
with $s<s^{\prime}.$ The continuous function $\alpha\left(  \cdot,x,x^{\prime
}\right)  $ restricted to $\partial X$ attains the negative value $\alpha
\bigl(g\left(  +\infty\right)  ,x,x^{\prime}\bigr)=-d\left(  x,x^{\prime
}\right)  $ and the positive value $\alpha\bigl(g\left(  -\infty\right)
,x,x^{\prime}\bigr)=d\left(  x,x^{\prime}\right)  .$ By connectivity of
$\partial X,$ $\exists$ a point $\xi\in\partial X$ such that $\alpha\left(
\xi,x,x^{\prime}\right)  =0.$

In order to apply theorem \ref{main}, we need to show that $\Omega=GY.$ For
this it suffices to show that the limit set $\Lambda_{X}\bigl(\Gamma\bigr)$ of
the action of $\Gamma$ on $X$ equals $\partial X$ (cf. proposition
\ref{omega}). $\Gamma$ acts on itself and the limit set $\Lambda_{\Gamma}%
\bigl(\Gamma\bigr)$ of this action is equal to $\partial\Gamma.$ Consider the
map $\Gamma\rightarrow X$ given by $\gamma\rightarrow\gamma\left(  p\right)  $
for some $p\in X$ fixed. This map is a quasi-isometry, hence induces a
homeomorphism $\partial\Gamma\rightarrow\partial X$ which takes $\Lambda
_{\Gamma}\bigl
(\Gamma\bigr)$ into $\Lambda_{X}\bigl(\Gamma\bigr)$. It follows that
$\Lambda_{X}\bigl(\Gamma\bigr)=\partial X$.\end{proof}

We proceed now to apply theorem \ref{main} to negatively curved polyhedra and
to $n$-dimen\-sional complete ideal polyhedra with curvature less
than or equal to $-1.$

\subsection{Negatively curved polyhedra}\label{negpolsec} 

A hyperbolic $n$-simplex is the convex hull, in hyperbolic $n$-space
$\mathbb{H}^{n},$ of $n+1$ points in general position. Let $Y$ be a locally
finite union of hyperbolic simplices glued together isometrically along faces
of the same dimension such that for every simplex $\sigma,$ each $\left(
k-1\right)  $-face of $\sigma$ is glued isometrically with some face
of some simplex of $Y.$

Such a space $Y$ is naturally a complete geodesic metric space (by results of
Bridson \cite{Bri} and Moussong \cite{Mou}) with distance function given as
follows: a broken geodesic from a point $x$ to a point $y$ is a map $f:\left[
a,b\right]  \rightarrow Y$ with $f\left(  a\right)  =x,$ $f\left(  b\right)
=y$ for which there exists a subdivision $a=t_{0}<t_{1}<\cdots<t_{k+1}=b$ of
$\left[  a,b\right]  $ such that for all $i=0,1,\ldots,k$ the restriction
$f|_{\,\left[  t_{i},t_{i+1}\right]  }$ is a geodesic whose image lies in a
single simplex. The length of a broken geodesic $f$ is defined to be
\[
\sum_{i=0}^{k}\mathfrak{\ell}\left(  f|_{\,\left[  t_{i},t_{i+1}\right]
}\right)  =\sum_{i=0}^{k}\left|  f\left(  t_{i}\right)  -f\left(
t_{i+1}\right)  \right|
\]
where the length inside a simplex is measured with respect to the hyperbolic
metric $\left|  \,\,\,\,\right|  .$ The distance $d\left(  x,y\right)  $ from
$x$ to $y$ is then defined to be the lower bound of the lengths of broken
geodesics from $x$ to $y.$

Since $Y$ is assumed to be locally finite, hence, locally compact, $Y$ becomes
a proper geodesic metric space.

\begin{definition}
\label{negpoldef}Such a space $Y$ is called a negatively curved polyhedron if
$Y$ with the induced length metric has curvature $\leq-1.$
\end{definition}

Recall that a geodesic metric space is said to have curvature less than or
equal to $\chi$ if each $y\in Y$ has a neighborhood $V_{y}$ such that every
geodesic triangle of perimeter strictly less than $\frac{2\pi}{\sqrt{\chi}}$
(=+$\infty$ when $\chi\leq0$) contained in $V_{y}$ satisfies $CAT\left(
\chi\right)  $.

Let $\tilde{Y}$ be the universal cover of $Y.$ Then $\tilde{Y}$ is a
$CAT\left(  -1\right)  $-space and $Y$ is the quotient $\tilde{Y}\!/\Gamma$
where $\Gamma$ is a discrete group of isometries of $\tilde{Y}$ isomorphic to
$\pi_{1}\left(  Y\right)  .$

\begin{notation}
\label{notation}It is explicit in the above definition that a negatively
curved polyhedron can be made up using simplices of various dimensions. Let
$\sigma$ be a $1$-dimen\-sional simplex in $Y$ such that $\sigma$ is not the
face of any $k$-simplex, $k\geq2 ,$ in $Y.$ Such a simplex will be called a
\textit{free }$1$-simplex. We will use in the sequel a subspace of the
$1$-skeleton of $Y$ which consists of all free $1$-simplices $\sigma$ in $Y.$
This subspace will be denoted by $Y^{\left[  1\right]  }$ and is not to be
confused with the $1$-skeleton of $Y.$ Observe that $Y^{\left[  1\right]  }$
may be empty. The (topological) boundary $Y^{\left[  1\right]  }%
\setminus\mathrm{Int\,}Y^{\left[  1\right]  }$ of $Y^{\left[  1\right]  },$
denoted by $bd\left(  Y^{\left[  1\right]  }\right)  ,$ is a discrete set of
points in $Y$ each of which is the $0$-face of some simplex of dimension
$k\geq2.$ 

For the universal cover $\tilde{Y}$ of $Y,$ the same
notation 
(i.e., $\tilde{Y}^{\left[  1\right]  })$ will be used.
\end{notation}

It is purely for convenience that we consider simplices of constant curvature
$-1$ instead of simplices of constant curvature $\chi,\chi<0.$ Moreover, we
may define our spaces to be negatively curved cell complexes. As any cell
complex can be made simplicial by subdivision, this involves no loss of
generality. For detailed definitions and properties of negatively curved
polyhedra we refer the reader to the treatments of Ballman \cite{B-G-H},
Bridson \cite{Bri} and Paulin \cite{Pau}.

\begin{theorem}
\label{mainncp}Let $Y$ be a negatively curved polyhedron which is not a graph.
Then the geodesic flow on $Y$ is topologically mixing, provided that the
non-wandering set $\Omega$ equals $GY$ and $\pi_{1}\left(  Y\right)  $ is
non-elementary. In particular, the geodesic flow on any compact negatively
curved polyhedron which is not a graph is topologically mixing.
\end{theorem}

As explained in section \ref{cex}, $1$-dimensional simplicial complexes
are of a special nature as far as topological mixing is concerned. This
continues to be the case with the zeros on the boundary of the generalized
Busemann function $\alpha.$ The next proposition, which asserts the existence
of such zeros, is false if the negatively curved polyhedron contains even a
single free $1$-simplex. This failure calls for specific treatment which is
given in the proof of theorem \ref{mainncp}.

\begin{proposition}
\label{alphazero}Let $Y$ be a negatively curved polyhedron and $\tilde{Y}$ its
universal cover. Then, $\forall x,y\in\tilde{Y}$ there exists $\xi\in
\partial\tilde{Y}$ such that $\alpha\left(  \xi,x,y\right)  =0$ provided that
the midpoint of the geodesic segment $\left[  x,y\right]  $ is not contained
in the interior of $\tilde{Y}^{\left[  1\right]  }$ (i.e., is not contained in
the interior of a $1$-simplex which is not the face of a $k$-simplex,
$k\geq2,$ of $\tilde{Y}).$
\end{proposition}

In the proof of the above proposition we will use the notion of the space of
directions: if $y$ is a point in $Y$ (or $\tilde{Y}$) we consider the space of
directions $D_{y}$ at the point $y.$ A point in $D_{y}$ is an equivalence
class of geodesic segments emanating from $y$ and angle measurement induces a
metric on $D_{y}.$ For details concerning angles in an arbitrary $CAT\left(
\chi\right)  $-space we refer the reader to \cite[Ch. I.3]{Bal}. In fact, the
simplicial structure of $Y$ induces a simplicial structure on $D_{y}$ so that
$D_{y}$ is locally a $CAT(1)$-space. If $x$ (resp. $\sigma$) is a point (resp.
a path) in $Y$ with $x\neq y$ (resp. $y\notin\mathrm{Im}\sigma$), then we
denote by $d\left(  x\right)  $ $\bigl
($resp. $d\left(  \sigma\right)  \bigr)$ the direction at $y$ pointing to $x$
(resp. the path $d\circ\sigma$ in $D_{y}$). For details concerning the space
of directions as well as for the following two facts needed in the sequel we
refer the reader to \cite[Chapter 10]{B-G-H}. Around any point $y$ there
exists a neighborhood $U_{y}$ such that%
\begin{equation}
\mathrm{If\,\,}\sigma\,\,\mathrm{is\,\,a\,\,geodesic\,\,segment\,\,in\,\,}%
U_{y}\mathrm{\,\,and\,\,}y\notin\mathrm{Im}\sigma,\mathrm{\,\,then\,\,}d\left(
\sigma\right)  \mathrm{\,\,has\,\,length\,\,}<\pi. \label{lysialpha}%
\end{equation}%
\begin{equation}\hspace{-1.5pc} 
\mathrm{If\,\,}\sigma\,\,\mathrm{is\,\,a\,\,geodesic\,\,segment\,\,in\,\,}%
U_{y},\mathrm{\,\,then\,\,}d\left(  \sigma\right)
\mathrm{\,\,is\,\,a\,\,geodesic\,\,segment\,\,in\,\,}D_{y}. \label{lysibita}%
\hspace{-1.5pc} \end{equation}

\begin{figure}[ptb]
\begin{center}
\font\thinlinefont=cmr5
\begingroup\makeatletter\ifx\SetFigFont\undefined%
\gdef\SetFigFont#1#2#3#4#5{%
  \reset@font\fontsize{#1}{#2pt}%
  \fontfamily{#3}\fontseries{#4}\fontshape{#5}%
  \selectfont}%
\fi\endgroup%
\mbox{\beginpicture
\setcoordinatesystem units <1.00000cm,1.00000cm>
\unitlength=1.00000cm
\linethickness=1pt
\setplotsymbol ({\makebox(0,0)[l]{\tencirc\symbol{'160}}})
\setshadesymbol ({\thinlinefont .})
\setlinear
\linethickness=1pt
\setplotsymbol ({\makebox(0,0)[l]{\tencirc\symbol{'160}}})
\putrule from  3.016 16.351 to  7.620 16.351
\plot  7.620 16.351  5.080 23.019 /
\plot  5.080 23.019  3.016 16.351 /
\linethickness=1pt
\setplotsymbol ({\makebox(0,0)[l]{\tencirc\symbol{'160}}})
\plot  5.080 23.019  3.651 25.718 /
\linethickness=1pt
\setplotsymbol ({\makebox(0,0)[l]{\tencirc\symbol{'160}}})
 \plot  3.016 16.351  3.55 23.971 /                                  
\linethickness=1pt
\setplotsymbol ({\makebox(0,0)[l]{\tencirc\symbol{'160}}})
 \plot  3.651 25.718  3.58 24.38 /
\linethickness= 0.500pt
\setplotsymbol ({\thinlinefont .})
\put{\makebox(0,0)[l]{\circle*{ 0.2}}} at  3.035 16.351
\linethickness= 0.500pt
\setplotsymbol ({\thinlinefont .})
\put{\makebox(0,0)[l]{\circle*{ 0.2}}} at  5.099 23.019
\linethickness= 0.500pt
\setplotsymbol ({\thinlinefont .})
\put{\makebox(0,0)[l]{\circle*{ 0.2}}} at  2.381 25.082
\linethickness= 0.500pt
\setplotsymbol ({\thinlinefont .})
\put{\makebox(0,0)[l]{\circle*{ 0.2}}} at  4.604 26.194
\linethickness= 0.500pt
\setplotsymbol ({\thinlinefont .})
\put{\makebox(0,0)[l]{\circle*{ 0.2}}} at  3.651 25.718
\linethickness= 0.500pt
\setplotsymbol ({\thinlinefont .})
\put{\makebox(0,0)[l]{\circle*{ 0.2}}} at  5.397 16.351
\linethickness= 0.500pt
\setplotsymbol ({\thinlinefont .})
\put{\makebox(0,0)[l]{\circle*{ 0.2}}} at  7.601 16.351
\linethickness=1pt
\setplotsymbol ({\makebox(0,0)[l]{\tencirc\symbol{'160}}})
\plot  2.381 25.082  4.604 26.194 /
\plot  4.604 26.194  5.080 23.019 /
\plot  5.080 23.019  2.381 25.082 /
\put{$x$} [lB] at  2.857 15.875
\put{$m$} [lB] at  5.397 15.875
\put{$y$} [lB] at  7.461 15.875
\put{$z_0$} [lB] at  5.556 23.019
\put{$y^\prime$} [lB] at  1.587 24.924
\put{$x^\prime$} [lB] at  5.080 26.194
\put{$z_1$} [lB] at  3.175 25.876
\linethickness=0pt
\putrectangle corners at  1.587 26.480 and  7.758 15.790
\endpicture}
\end{center}
\caption{}%
\label{monotonicity}%
\end{figure}

\begin{proof}[Proof of Proposition \ref{alphazero}] Let $x,y\in\tilde{Y}$
with $x\neq y$ and denote by $m$ the midpoint of the geodesic segment $\left[
x,y\right]  $. Consider the set
\[
Z=\left\{  z\in\tilde{Y}\bigm\vert\alpha\left(  z,x,y\right)  =0\right\}.
\]\enlargethispage{\baselineskip}
Observe that $Z\neq\emptyset,$ as $m\in Z.$ It suffices to show that $Z$ is
not bounded in $\tilde{Y}.$ For, if $\left\{  z_{n}\right\}  _{n\in\mathbb{N}%
}$ is a sequence in $Z$ with $d\left(  z_{n},m\right)  \rightarrow+\infty$ as
$n\rightarrow\infty$, then by choosing, if necessary, a subsequence we have
that $\left\{  z_{n}\right\}  _{n\in\mathbb{N}}$ converges to some point
$\xi\in\partial\tilde{Y}$ and by continuity of the $\alpha$ function,
\[
\alpha\left(  \xi,x,y\right)  =\mathrm{lim}_{t\rightarrow\infty}\alpha\bigl
(z_{n},x,y\bigr)=0.
\]
We proceed to show that $Z$ is not bounded in $\tilde{Y}.$ Assume, on the
contrary, that there exists a point $z_{0}\in Z$ of maximal distance from $x$
(and, hence, from $y$). We have two cases: \vspace{2mm}

\noindent\textit{Case} A: $z_{0}=m.$
\vspace{2mm}

Consider
the space of directions $D_{m}$ at the point $m.$ If $D_{m}$ is connected, then
for points $m_{1},m_{2}$ on the geodesic segments\pagebreak\ $\left[  m,x\right]  $ and
$\left[  m,y\right]  ,$ respectively, which are arbitrarily close to $m,$
there exists a path $\tau$ in $D_{m}$ with endpoints $d\left(  m_{1}\right)  $
and $d\left(  m_{2}\right)  $. This path $\tau$ determines a path $\sigma$ in
$\tilde{Y}$ with endpoints $m_{1}$ and $m_{2}$ so that $\tau=d\left(
\sigma\right)  $ and $\mathrm{Im}\sigma$ does not contain $m.$ The generalized
Busemann function $\alpha$ attains negative and positive values on the
endpoints of $\sigma.$ Thus, there exists a point $z_{1}$ on $\mathrm{Im}%
\sigma$ with $\alpha\left(  z_{1},x,y\right)  =0.$ The piecewise geodesic
$\left[  x,z_{1}\right]  \cup\left[  z_{1},y\right]  $ has length, by
uniqueness of (length minimizing) geodesics in $\tilde{Y},$ strictly bigger
than $d\left(  x,y\right)  =2d\left(  x,m\right)  =2d\left(  y,m\right)  .$
Thus,
\[
d\left(  z_{1},x\right)  =d\left(  z_{1},y\right)  \gneqq d\left(
x,m\right)
\]
which contradicts the fact that $m=z_{0}$ is a point in $Z$ of maximal
distance from $x.$\newline \noindent Assume now that $D_{m}$ is not connected.
We only have to deal with the case in which the geodesic segments $\left[
m,x\right] $ and $\left[  m,y\right]  $ determine points $d\left(  x\right)  $
and $d\left(  y\right)  $ in distinct components, say $C_{x}$ and $C_{y},$ of
$D_{m}.$ If there exists a connected component $C_{0}$ of $D_{m}$ distinct
from $C_{x}$ and $C_{y}$, then pick a point $z_{1}$ so that $d\left(
z_{1}\right)  \in C_{0}.$ Then the geodesic segments $\left[  x,z_{1}\right]
$ and $\left[  y,z_{1}\right]  $ necessarily contain $z_{0}.$ It follows that
$\alpha\left(  z_{1},x,y\right)  =0$ and
\[
d\left(  z_{1},x\right) = d\left(  z_{1},m\right) + d\left( m,x\right) >
    d\left( m,x\right)
\]
which contradicts the fact that $m=z_{0}$ is a point in $Z$
of maximal distance from $x.$ \newline Assume now that $C_{x}$ and $C_{y}$ are
the only components of $D_{m}.$ By the midpoint assumption on $\left[
x,y\right]  ,$ the sets $C_{x}$ and $C_{y}$ cannot be both singletons. Assume
$C_{y}$ is not a singleton. $C_{y}$ satisfies $CAT(1)$-inequality and,
therefore, the \textit{systole} of $Link(m,\widetilde{Y})$ is greater than or equal
to $2\pi$ (see \cite[theorem 3.15]{Pau}). Hence, we can extend $\left[
y,m\right]  $ to a geodesic segment $\left[  y,z_{1}\right]  $ which contains
$m$ in its interior so that $d\left(  z_{1}\right)  \in C_{y}.$ Then $\left[
x,m\right]  \cup\left[  m,z_{1}\right]  $ and $\left[  y,m\right]  \cup\left[
m,z_{1}\right]  $ are both geodesic segments. It now follows that
\[
d\left(  z_{1},x\right)  =d\left(  z_{1},y\right)  \gneqq d\left(
m,x\right)
\]
which completes the proof in this case.\vspace{2mm}\newline \noindent
\textit{Case }B: $z_{0}\neq m.$

Choose a
neighborhood $U$ around $z_{0}$ so that statements (\ref{lysialpha}) and
(\ref{lysibita}) hold. Extend the geodesic segments $\left[  x,z_{0}\right]  $
and $\left[  y,z_{0}\right]  ,$ i.e., choose points $x^{\prime},y^{\prime}\in
U\setminus\left\{  z_{0}\right\}  $ so that $\left[  x,z_{0}\right]
\cup\left[  z_{0},x^{\prime}\right]  $ and $\left[  y,z_{0}\right]
\cup\left[  z_{0},y^{\prime}\right]  $ are both geodesic segments. Amongst all
possible choices for the pair $x^{\prime}$ and $y^{\prime}$ pick one so that
\[
\measuredangle\,_{z_{0}}\left(  x^{\prime},y^{\prime}\right)  <\pi.
\]
If $d\left(  x^{\prime}\right)  =d\left(  y^{\prime}\right)  $, then for some
point $z_{1}$ on $\left[  z_{0},x^{\prime}\right]  \cap\left[  z_{0}%
,y^{\prime}\right]  \setminus\left\{  z_{0}\right\}  $ we have
\[
d\left(  z_{1},x\right)  =d\left(  z_{1},y\right)  \gvertneqq d\left(
z_{0},x\right)
\]
which implies that $z_{0}$ is not of maximal distance from $x.$ Hence, we may
assume that $d\left(  x^{\prime}\right)  \neq d\left(  y^{\prime}\right)
$.

It is clear that $\alpha\left(  x^{\prime},x,y\right)  \leq0$ and
$\alpha\left(  y^{\prime},x,y\right)  \geq0$. If either $\alpha\left(
x^{\prime},x,y\right)  =0$ or $\alpha\left(  y^{\prime},x,y\right)  =0$, then,
again, the point $z_{0}$ is not of maximal distance from $x.$ Thus we may
assume that $\alpha\left(  x^{\prime},x,y\right)  <0$ and $\alpha\left(
y^{\prime},x,y\right)  >0$ which implies that
\[
\exists\,\,z_{1}\in\left[  x^{\prime},y^{\prime}\right]  :\alpha\left(
z_{1},x,y\right)  =0.
\]
We proceed to show that $d\left(  z_{1},x\right)  >d\left(  z_{0},x\right)  .$
Denote by $\sigma,\sigma_{1}$ and $\sigma_{2}$ the geodesic segments $\left[
x^{\prime},y^{\prime}\right]  ,\left[  x^{\prime},z_{1}\right]  $ and $\left[
z_{1},y^{\prime}\right]  $ respectively. Observe that $z_{0}\notin
\mathrm{Im}\sigma$ because $x^{\prime},y^{\prime}$ are chosen so that
$\measuredangle\,_{z_{0}}\left(  x^{\prime},y^{\prime}\right)  <\pi.$ By
(\ref{lysibita}) the projection $d\left(  \sigma\right)  $ of $\sigma$ in the
space of directions $D_{z_{0}}$ is a geodesic segment, hence,
\[
d\left(  \sigma_{1}\right)  +d\left(  \sigma_{2}\right)  =d\left(
\sigma\right)  .
\]
As the length of $d\left(  \sigma\right)  $ is smaller, then $\pi,$ either
$d\left(  \sigma_{1}\right)  $ or, $d\left(  \sigma_{2}\right)  $ has length
$<\pi/2.$ In other words, either $\measuredangle\,_{z_{0}}\left(
z_{1},x^{\prime}\right)  <\pi/2$ or, $\measuredangle\,_{z_{0}}\left(
z_{1},y^{\prime}\right)  <\pi/2.$ We may assume that
\begin{equation}
\measuredangle\,_{z_{0}}\left(  z_{1},x^{\prime}\right)  <\pi/2 \label{star}
\end{equation}
(if $\measuredangle\,_{z_{0}}\left(  z_{1},y^{\prime}\right)  <\pi/2$ we
proceed in an identical way). As $\left[  x,z_{0}\right]  \cup\left[
z_{0},x^{\prime}\right]  $ is a geodesic segment, $\measuredangle\,_{z_{0}%
}\left(  x,x^{\prime}\right)  \geq\pi.$ Since
\[
\measuredangle\,_{z_{0}}\left(  x,z_{1}\right)  +\measuredangle\,_{z_{0}%
}\left(  z_{1},x^{\prime}\right)  \geq\measuredangle\,_{z_{0}}\left(
x,x^{\prime}\right)
\]
we have, by (\ref{star}), that
\[
\measuredangle\,_{z_{0}}\left(  x,z_{1}\right)  >\pi/2.
\]
The latter inequality implies that $d\left(  z_{1},x\right)  >d\left(
z_{0},x\right)  $ which completes the proof of the proposition.\end{proof}

In the proof of proposition \ref{mainncp} we will need the following
construction:

\begin{construc}
Let $Y$ be a negatively curved polyhedron.
Consider the subspace $Y^{\left[  1\right]  }$ of $Y$ (explained in
\ref{notation} above) and let $Y_{i},i\in\mathbb{N}$ be the connected
components of $Y\setminus\mathrm{Int\,}Y^{\left[  1\right]  }.$ For each $i,$
let $\left\{  \sigma_{i_{j}}\bigm\vert
j=1,2,\ldots\right\}  $ be an enumeration (possibly infinite) of the set
$Y_{i}\cap Y^{\left[  1\right]  }.$ We glue the components $Y_{i}%
,i\in\mathbb{N}$ together according to the following rule:

\begin{quotation}
The $0$-face $\sigma_{i_{j}}$ of $Y_{i}$ is identified with the $0$-face
$\sigma_{i_{j^{\prime}}^{\prime}}$ of $Y_{i^{\prime}}$ if there exists a path
lying entirely in $Y^{\left[  1\right]  }$ with endpoints $\sigma_{i_{j}} $
and $\sigma_{i_{j^{\prime}}^{\prime}}.$
\end{quotation}

Note that the equality $i=i^{\prime}$ is allowed in the above rule. In this
way we obtain a negatively curved polyhedron denoted by $Y_{-1}$ which does
not contain any $1$-dimensional simplices. The image of a geodesic $g$ in $GY$
determines a unique geodesic line in $Y_{-1}.$ By employing base points on $Y$
and $Y_{-1}$ each geodesic $g$ in $GY$ determines a unique geodesic $g_{-1}$
in $GY_{-1}.$ This map
\begin{equation}
GY\rightarrow GY_{-1}\mathrm{\,\,is\,\,surjective.} \label{minusonesurj}%
\end{equation}
Given a subset $\mathcal{V}$ (resp. a point $f$) in $GY$ we will be denoting
by $\mathcal{V}_{-1}$ (resp. $f_{-1}$) the corresponding subset (resp. point)
in $GY_{-1}$ under the above map.
\end{construc}

\begin{proof}[Proof of Theorem \ref{mainncp}] First observe that in the
case of a compact negatively curved polyhedron the assumption $\Omega=GY$
follows exactly as in the proof of Corollary \ref{Bestvina} above. Moreover,
as $\pi_{1}\left(  Y\right)  $ acts co-compactly on the hyperbolic space
$\tilde{Y}$ it follows (by a theorem of Gromov) that $\pi_{1}\left(  Y\right)
$ is a hyperbolic group, hence, non-elementary. Therefore, theorem
\ref{mainncp} can be stated for compact negatively curved polyhedra without
any hypothesis at all.

If $Y$ does not contain $1$-dimensional simplices, then the conclusion of the
theorem follows from theorem \ref{main} and proposition \ref{alphazero}. In
order to deal with the general case we will modify the proof of theorem
\ref{main} at the point where assumption (1) of theorem \ref{main} is used.
This modification will assert that the midpoint assumption of proposition
\ref{alphazero} is fulfilled, hence, the zeros for the $\alpha$ function
needed do, in fact, exist.

Recall that assumption (1) of theorem \ref{main} is only used in the proof of
proposition \ref{ss} where given arbitrary open sets $\mathcal{O}%
,\mathcal{U}\subseteq GY$ we choose $f\in p^{-1}\left(  \mathcal{O}\right)  $
and $h\in p^{-1}\left(  \mathcal{U}\right)  $ such that
\[
\exists \,\, a\,\, hyperbolic\,\, \phi \in\Gamma :
\bigl(f\left(  +\infty\right)  ,h\left(  +\infty\right)  \bigr)=\left(
\phi\left(  +\infty\right)  ,\phi\left(  -\infty\right)  \right)
\]
and then assumption (1) is employed to obtain, for each $n,$ a point $\xi_{n}
$ in $\partial X$ such that
\[
\alpha\bigl(\xi_{n},f\left(  0\right)  ,\phi^{n}\left(  h\left(  0\right)
\right)  \bigr)=0.
\]
We only need to show that by an appropriate choice of $f,h$ and $\phi$ the
midpoint of the geodesic segment $\left[  f\left(  0\right)  ,\phi^{n}\left(
h\left(  0\right)  \right)  \right]  $ lies, for all $n$ large enough, in the
interior of a $k$-simplex of $\tilde{Y}$ with $k\geq2.$ Let $\mathcal{O}%
,\mathcal{U}\subseteq GY$ be given. We may choose a geodesic $h^{\prime}%
\in\mathcal{U},$ a compact set $K_{h^{\prime}}\subset\mathbb{R}$ and a real
$\varepsilon_{h^{\prime}}>0$ such that

\begin{description} 
\item [(a)]the neighborhood $\mathcal{U}^{\prime}$ around $h^{\prime}$
determined by $K_{h^{\prime}}$ and $\varepsilon_{h^{\prime}},$ i.e.
\newline \hspace*{4mm}$\mathcal{U}^{\prime}=\left\{  g\in GY\bigm\vert d\bigl
(g\left(  t\right)  ,h^{\prime}\left(  t\right)  \bigr)<\varepsilon
_{h^{\prime}}\,\,\forall\,\,t\in K_{h^{\prime}}\right\}  ,$ is a subset of
$\mathcal{U}.$
\end{description}
We may refine the choices of $h^{\prime},K_{h^{\prime}}$ and
$\varepsilon_{h^{\prime}}$ so that, in addition, the following property is
satisfied:

\begin{description}
\item [(b)]for all $t>\max K_{h^{\prime}},$ $h^{\prime}\left(  t\right)  $
lies in a single component, say $Y_{0},$ \newline \hspace*{4mm}of
$Y\setminus\mathrm{Int\,}Y^{\left[  1\right]  }$
\end{description}
as follows: we need to consider the time $t_{h^{\prime}}\notin
K_{h^{\prime}}$ at which the geodesic $h^{\prime}$ first enters the set
$Y^{\left[  1\right]  },$ i.e., set
\[
t_{h^{\prime}}=\mathrm{inf}\left\{  t\in\left[  \max K_{h^{\prime}}%
,+\infty\right) \!\! \bigm\vert\!\! \,h^{\prime}\left(  t\right)  \in Y^{\left[
1\right]  }\,\,\mathrm{and\,\,}h^{\prime}\left(  t+\varepsilon\right)
\in\mathrm{Int\,}Y^{\left[  1\right]  }\,\,\forall\varepsilon
>0\,\,\mathrm{small}\right\}\!\!.
\]
If $t_{h^{\prime}}=+\infty,$ i.e., $h^{\prime}\bigl(\left[  \max K_{h^{\prime
}},+\infty\right)  \bigr)$ does not intersect $\mathrm{Int\,}Y^{\left[
1\right]  },$ then $h^{\prime}|_{\left[  \max K_{h^{\prime}},+\infty\right)
}$ stays in a single component, say $Y_{0},$ of $Y\setminus\mathrm{Int\,}%
Y^{\left[  1\right]  }$ and, hence, property \textbf{(b)} is satisfied.
Suppose now that $t_{h^{\prime}}\neq+\infty.$ Extend (in an arbitrary way) the
geodesic ray $h^{\prime}:\left(  -\infty,t_{h^{\prime}}\right]  \rightarrow Y$
to a geodesic ray $\left(  -\infty,T\right]  \rightarrow Y,$ for some
$T>t_{h^{\prime}},$ denoted again by $h^{\prime},$ so that $h^{\prime}\left(
T\right)  \in bd\left(  Y^{\left[  1\right]  }\right)  .$ Let $Y_{0}$ be the
component of $Y\setminus\mathrm{Int\,}Y^{\left[  1\right]  }$ which contains
$h^{\prime}\left(  T\right)  .$ We may now extend the geodesic ray $\left(
-\infty,T\right]  \rightarrow Y$ to a geodesic line, denoted again by
$h^{\prime},$ so that
\[
h^{\prime}\bigl(\left[  T,+\infty\right)  \bigr)\subset Y_{0}.
\]
By enlarging, if necessary, the compact set $K_{h^{\prime}}$ to contain $T$
the choice of $h^{\prime}\in\mathcal{U},$ $K_{h^{\prime}}\subset\mathbb{R}$
and $\varepsilon_{h^{\prime}}>0$ satisfying \textbf{(a)} and \textbf{(b)} is
complete. We need to do the same thing for $\mathcal{O},$ i.e., to choose a
geodesic $f^{\prime}\in\mathcal{O},$ a compact set $K_{f^{\prime}}%
\subset\mathbb{R}$ and a real $\varepsilon_{f^{\prime}}>0$ such that

\begin{description}
\item [(c)]the neighborhood $\mathcal{O}^{\prime}$ around $f^{\prime}$
determined by $K_{f^{\prime}}$ and $\varepsilon_{f^{\prime}}$ \newline
\hspace*{4mm}is a subset of $\mathcal{O},$ 

\item[(d)] for all $t>\max K_{f^{\prime}},$ $f^{\prime}\left(  t\right)  $
lies in the same component $Y_{0}$ of \newline \hspace*{4mm}$Y\setminus
\mathrm{Int\,}Y^{\left[  1\right]  }.$
\end{description}

\begin{figure}[ptb]
\begin{center}
\font\thinlinefont=cmr5
\begingroup\makeatletter\ifx\SetFigFont\undefined%
\gdef\SetFigFont#1#2#3#4#5{%
  \reset@font\fontsize{#1}{#2pt}%
  \fontfamily{#3}\fontseries{#4}\fontshape{#5}%
  \selectfont}%
\fi\endgroup%
\mbox{\beginpicture
\setcoordinatesystem units <1.00000cm,1.00000cm>
\unitlength=1.00000cm
\linethickness=1pt
\setplotsymbol ({\makebox(0,0)[l]{\tencirc\symbol{'160}}})
\setshadesymbol ({\thinlinefont .})
\setlinear
\linethickness= 0.500pt
\setplotsymbol ({\thinlinefont .})
\put{\makebox(0,0)[l]{\circle*{ 0.2}}} at  3.670 24.448
\linethickness= 0.500pt
\setplotsymbol ({\thinlinefont .})
\put{\makebox(0,0)[l]{\circle*{ 0.2}}} at  3.670 21.749
\linethickness= 0.500pt
\setplotsymbol ({\thinlinefont .})
\put{\makebox(0,0)[l]{\circle*{ 0.2}}} at  7.480 18.574
\linethickness= 0.500pt
\setplotsymbol ({\thinlinefont .})
\put{\makebox(0,0)[l]{\circle*{ 0.2}}} at  7.461 21.749
\linethickness= 0.500pt
\setplotsymbol ({\thinlinefont .})
\put{\makebox(0,0)[l]{\circle*{ 0.2}}} at  7.461 20.320
\linethickness= 0.500pt
\setplotsymbol ({\thinlinefont .})
\put{\makebox(0,0)[l]{\circle*{ 0.2}}} at  5.893 24.765
\linethickness= 0.500pt
\setplotsymbol ({\thinlinefont .})
\put{\makebox(0,0)[l]{\circle*{ 0.2}}} at  9.544 21.749
\linethickness= 0.500pt
\setplotsymbol ({\thinlinefont .})
\put{\makebox(0,0)[l]{\circle*{ 0.2}}} at  3.632 23.178
\linethickness= 0.500pt
\setplotsymbol ({\thinlinefont .})
\put{\makebox(0,0)[l]{\circle*{ 0.2}}} at  9.506 23.495
\linethickness= 0.500pt
\setplotsymbol ({\thinlinefont .})
\putrule from 13.335 21.749 to 13.176 21.749
\linethickness= 0.500pt
\setplotsymbol ({\thinlinefont .})
\putrule from  7.461 18.500 to  7.461 21.749
\linethickness=1pt
\setplotsymbol ({\makebox(0,0)[l]{\tencirc\symbol{'160}}})
\putrule from  1.746 21.749 to 13.335 21.749
\linethickness=1pt
\setplotsymbol ({\makebox(0,0)[l]{\tencirc\symbol{'160}}})
\putrule from  9.525 23.495 to  9.525 21.749
\linethickness=1pt
\setplotsymbol ({\makebox(0,0)[l]{\tencirc\symbol{'160}}})
\putrule from  3.651 24.448 to  3.651 21.749
\linethickness=1pt
\setplotsymbol ({\makebox(0,0)[l]{\tencirc\symbol{'160}}})
\plot  3.624 23.139  5.370 21.869 /
\linethickness=1pt
\setplotsymbol ({\makebox(0,0)[l]{\tencirc\symbol{'160}}})
\plot  5.605 21.656  7.351 20.386 /
\linethickness=1pt
\setplotsymbol ({\makebox(0,0)[l]{\tencirc\symbol{'160}}})
\setdashes < 0.2858cm>
\plot  1.746 21.431  1.750 21.431 /
\plot  1.750 21.431  1.763 21.429 /
\plot  1.763 21.429  1.782 21.427 /
\plot  1.782 21.427  1.812 21.423 /
\plot  1.812 21.423  1.852 21.416 /
\plot  1.852 21.416  1.903 21.410 /
\plot  1.903 21.410  1.960 21.402 /
\plot  1.960 21.402  2.026 21.391 /
\plot  2.026 21.391  2.095 21.380 /
\plot  2.095 21.380  2.167 21.370 /
\plot  2.167 21.370  2.237 21.359 /
\plot  2.237 21.359  2.307 21.349 /
\plot  2.307 21.349  2.373 21.338 /
\plot  2.373 21.338  2.436 21.328 /
\plot  2.436 21.328  2.496 21.317 /
\plot  2.496 21.317  2.551 21.306 /
\plot  2.551 21.306  2.601 21.296 /
\plot  2.601 21.296  2.650 21.283 /
\plot  2.650 21.283  2.699 21.273 /
\plot  2.699 21.273  2.745 21.260 /
\plot  2.745 21.260  2.794 21.247 /
\plot  2.794 21.247  2.841 21.234 /
\plot  2.841 21.234  2.889 21.220 /
\plot  2.889 21.220  2.936 21.205 /
\plot  2.936 21.205  2.982 21.190 /
\plot  2.982 21.190  3.029 21.175 /
\plot  3.029 21.175  3.076 21.160 /
\plot  3.076 21.160  3.122 21.143 /
\plot  3.122 21.143  3.167 21.129 /
\plot  3.167 21.129  3.213 21.112 /
\plot  3.213 21.112  3.260 21.097 /
\plot  3.260 21.097  3.308 21.080 /
\plot  3.308 21.080  3.355 21.063 /
\plot  3.355 21.063  3.404 21.046 /
\plot  3.404 21.046  3.450 21.029 /
\plot  3.450 21.029  3.501 21.012 /
\plot  3.501 21.012  3.550 20.993 /
\plot  3.550 20.993  3.600 20.974 /
\plot  3.600 20.974  3.651 20.955 /
\plot  3.651 20.955  3.698 20.936 /
\plot  3.698 20.936  3.747 20.917 /
\plot  3.747 20.917  3.795 20.898 /
\plot  3.795 20.898  3.844 20.879 /
\plot  3.844 20.879  3.893 20.860 /
\plot  3.893 20.860  3.941 20.841 /
\plot  3.941 20.841  3.992 20.820 /
\plot  3.992 20.820  4.041 20.800 /
\plot  4.041 20.800  4.092 20.779 /
\plot  4.092 20.779  4.142 20.760 /
\plot  4.142 20.760  4.191 20.739 /
\plot  4.191 20.739  4.242 20.718 /
\plot  4.242 20.718  4.293 20.697 /
\plot  4.293 20.697  4.343 20.676 /
\plot  4.343 20.676  4.396 20.654 /
\plot  4.396 20.654  4.447 20.631 /
\plot  4.447 20.631  4.498 20.608 /
\plot  4.498 20.608  4.551 20.585 /
\plot  4.551 20.585  4.604 20.559 /
\plot  4.604 20.559  4.655 20.534 /
\plot  4.655 20.534  4.710 20.506 /
\plot  4.710 20.506  4.763 20.479 /
\plot  4.763 20.479  4.809 20.453 /
\plot  4.809 20.453  4.856 20.428 /
\plot  4.856 20.428  4.902 20.403 /
\plot  4.902 20.403  4.949 20.375 /
\plot  4.949 20.375  4.995 20.348 /
\plot  4.995 20.348  5.044 20.318 /
\plot  5.044 20.318  5.093 20.290 /
\plot  5.093 20.290  5.139 20.261 /
\plot  5.139 20.261  5.188 20.231 /
\plot  5.188 20.231  5.237 20.201 /
\plot  5.237 20.201  5.285 20.172 /
\plot  5.285 20.172  5.334 20.142 /
\plot  5.334 20.142  5.383 20.113 /
\plot  5.383 20.113  5.431 20.081 /
\plot  5.431 20.081  5.480 20.049 /
\plot  5.480 20.049  5.529 20.019 /
\plot  5.529 20.019  5.580 19.988 /
\plot  5.580 19.988  5.628 19.956 /
\plot  5.628 19.956  5.679 19.924 /
\plot  5.679 19.924  5.728 19.890 /
\plot  5.728 19.890  5.779 19.859 /
\plot  5.779 19.859  5.827 19.825 /
\plot  5.827 19.825  5.878 19.791 /
\plot  5.878 19.791  5.929 19.757 /
\plot  5.929 19.757  5.982 19.721 /
\plot  5.982 19.721  6.032 19.685 /
\plot  6.032 19.685  6.081 19.651 /
\plot  6.081 19.651  6.130 19.615 /
\plot  6.130 19.615  6.179 19.579 /
\plot  6.179 19.579  6.229 19.543 /
\plot  6.229 19.543  6.280 19.505 /
\plot  6.280 19.505  6.331 19.467 /
\plot  6.331 19.467  6.384 19.429 /
\plot  6.384 19.429  6.437 19.389 /
\plot  6.437 19.389  6.490 19.348 /
\plot  6.490 19.348  6.545 19.308 /
\plot  6.545 19.308  6.598 19.268 /
\plot  6.598 19.268  6.653 19.228 /
\plot  6.653 19.228  6.708 19.185 /
\plot  6.708 19.185  6.763 19.145 /
\plot  6.763 19.145  6.816 19.103 /
\plot  6.816 19.103  6.871 19.061 /
\plot  6.871 19.061  6.924 19.020 /
\plot  6.924 19.020  6.979 18.978 /
\plot  6.979 18.978  7.032 18.938 /
\plot  7.032 18.938  7.082 18.895 /
\plot  7.082 18.895  7.133 18.855 /
\plot  7.133 18.855  7.184 18.815 /
\plot  7.184 18.815  7.233 18.773 /
\plot  7.233 18.773  7.281 18.733 /
\plot  7.281 18.733  7.328 18.692 /
\plot  7.328 18.692  7.374 18.654 /
\plot  7.374 18.654  7.419 18.614 /
\plot  7.419 18.614  7.461 18.574 /
\plot  7.461 18.574  7.506 18.531 /
\plot  7.506 18.531  7.550 18.487 /
\plot  7.550 18.487  7.595 18.443 /
\plot  7.595 18.443  7.639 18.396 /
\plot  7.639 18.396  7.686 18.347 /
\plot  7.686 18.347  7.732 18.294 /
\plot  7.732 18.294  7.779 18.239 /
\plot  7.779 18.239  7.830 18.182 /
\plot  7.830 18.182  7.880 18.123 /
\plot  7.880 18.123  7.931 18.062 /
\plot  7.931 18.062  7.984 17.998 /
\plot  7.984 17.998  8.037 17.935 /
\plot  8.037 17.935  8.088 17.869 /
\plot  8.088 17.869  8.139 17.808 /
\plot  8.139 17.808  8.187 17.748 /
\plot  8.187 17.748  8.232 17.691 /
\plot  8.232 17.691  8.274 17.640 /
\plot  8.274 17.640  8.310 17.594 /
\plot  8.310 17.594  8.340 17.556 /
\plot  8.340 17.556  8.365 17.524 /
\plot  8.365 17.524  8.384 17.501 /
\plot  8.384 17.501  8.397 17.484 /
\plot  8.397 17.484  8.407 17.471 /
\plot  8.407 17.471  8.412 17.465 /
\plot  8.412 17.465  8.414 17.462 /
\linethickness=1pt
\setplotsymbol ({\makebox(0,0)[l]{\tencirc\symbol{'160}}})
\plot 13.310 22.032 13.305 22.032 /
\plot 13.305 22.032 13.293 22.035 /
\plot 13.293 22.035 13.274 22.037 /
\plot 13.274 22.037 13.244 22.041 /
\plot 13.244 22.041 13.204 22.047 /
\plot 13.204 22.047 13.153 22.054 /
\plot 13.153 22.054 13.096 22.062 /
\plot 13.096 22.062 13.030 22.073 /
\plot 13.030 22.073 12.960 22.083 /
\plot 12.960 22.083 12.888 22.094 /
\plot 12.888 22.094 12.819 22.104 /
\plot 12.819 22.104 12.749 22.115 /
\plot 12.749 22.115 12.683 22.126 /
\plot 12.683 22.126 12.620 22.136 /
\plot 12.620 22.136 12.560 22.147 /
\plot 12.560 22.147 12.505 22.157 /
\plot 12.505 22.157 12.454 22.168 /
\plot 12.454 22.168 12.406 22.181 /
\plot 12.406 22.181 12.357 22.191 /
\plot 12.357 22.191 12.311 22.204 /
\plot 12.311 22.204 12.262 22.217 /
\plot 12.262 22.217 12.215 22.229 /
\plot 12.215 22.229 12.167 22.244 /
\plot 12.167 22.244 12.120 22.259 /
\plot 12.120 22.259 12.073 22.274 /
\plot 12.073 22.274 12.027 22.288 /
\plot 12.027 22.288 11.980 22.303 /
\plot 11.980 22.303 11.934 22.320 /
\plot 11.934 22.320 11.889 22.335 /
\plot 11.889 22.335 11.843 22.352 /
\plot 11.843 22.352 11.796 22.367 /
\plot 11.796 22.367 11.748 22.384 /
\plot 11.748 22.384 11.701 22.401 /
\plot 11.701 22.401 11.652 22.418 /
\plot 11.652 22.418 11.606 22.435 /
\plot 11.606 22.435 11.555 22.451 /
\plot 11.555 22.451 11.506 22.471 /
\plot 11.506 22.471 11.455 22.490 /
\plot 11.455 22.490 11.405 22.509 /
\plot 11.405 22.509 11.358 22.528 /
\plot 11.358 22.528 11.309 22.547 /
\plot 11.309 22.547 11.261 22.566 /
\plot 11.261 22.566 11.212 22.585 /
\plot 11.212 22.585 11.163 22.604 /
\plot 11.163 22.604 11.115 22.623 /
\plot 11.115 22.623 11.064 22.644 /
\plot 11.064 22.644 11.015 22.663 /
\plot 11.015 22.663 10.964 22.684 /
\plot 10.964 22.684 10.914 22.703 /
\plot 10.914 22.703 10.865 22.725 /
\plot 10.865 22.725 10.814 22.746 /
\plot 10.814 22.746 10.763 22.767 /
\plot 10.763 22.767 10.712 22.788 /
\plot 10.712 22.788 10.660 22.809 /
\plot 10.660 22.809 10.609 22.832 /
\plot 10.609 22.832 10.558 22.856 /
\plot 10.558 22.856 10.505 22.879 /
\plot 10.505 22.879 10.452 22.904 /
\plot 10.452 22.904 10.401 22.930 /
\plot 10.401 22.930 10.346 22.957 /
\plot 10.346 22.957 10.293 22.985 /
\plot 10.293 22.985 10.247 23.010 /
\plot 10.247 23.010 10.200 23.036 /
\plot 10.200 23.036 10.154 23.061 /
\plot 10.154 23.061 10.107 23.089 /
\plot 10.107 23.089 10.061 23.116 /
\plot 10.061 23.116 10.012 23.146 /
\plot 10.012 23.146  9.963 23.173 /
\plot  9.963 23.173  9.917 23.203 /
\plot  9.917 23.203  9.868 23.233 /
\plot  9.868 23.233  9.819 23.262 /
\plot  9.819 23.262  9.771 23.292 /
\plot  9.771 23.292  9.722 23.321 /
\plot  9.722 23.321  9.673 23.351 /
\plot  9.673 23.351  9.624 23.383 /
\plot  9.624 23.383  9.576 23.415 /
\plot  9.576 23.415  9.527 23.444 /
\plot  9.527 23.444  9.476 23.476 /
\plot  9.476 23.476  9.428 23.508 /
\plot  9.428 23.508  9.377 23.539 /
\plot  9.377 23.539  9.328 23.573 /
\plot  9.328 23.573  9.277 23.605 /
\plot  9.277 23.605  9.229 23.639 /
\plot  9.229 23.639  9.178 23.673 /
\plot  9.178 23.673  9.127 23.707 /
\plot  9.127 23.707  9.074 23.743 /
\plot  9.074 23.743  9.023 23.779 /
\plot  9.023 23.779  8.975 23.812 /
\plot  8.975 23.812  8.926 23.848 /
\plot  8.926 23.848  8.877 23.884 /
\plot  8.877 23.884  8.826 23.920 /
\plot  8.826 23.920  8.776 23.959 /
\plot  8.776 23.959  8.723 23.997 /
\plot  8.723 23.997  8.672 24.037 /
\plot  8.672 24.037  8.619 24.077 /
\plot  8.619 24.077  8.566 24.117 /
\plot  8.566 24.117  8.511 24.158 /
\plot  8.511 24.158  8.456 24.200 /
\plot  8.456 24.200  8.403 24.242 /
\plot  8.403 24.242  8.348 24.282 /
\plot  8.348 24.282  8.293 24.325 /
\plot  8.293 24.325  8.238 24.367 /
\plot  8.238 24.367  8.183 24.409 /
\plot  8.183 24.409  8.130 24.452 /
\plot  8.130 24.452  8.077 24.494 /
\plot  8.077 24.494  8.024 24.534 /
\plot  8.024 24.534  7.971 24.577 /
\plot  7.971 24.577  7.921 24.617 /
\plot  7.921 24.617  7.870 24.657 /
\plot  7.870 24.657  7.821 24.697 /
\plot  7.821 24.697  7.772 24.737 /
\plot  7.772 24.737  7.726 24.776 /
\plot  7.726 24.776  7.681 24.814 /
\plot  7.681 24.814  7.637 24.852 /
\plot  7.637 24.852  7.595 24.890 /
\plot  7.595 24.890  7.546 24.934 /
\plot  7.546 24.934  7.499 24.977 /
\plot  7.499 24.977  7.451 25.023 /
\plot  7.451 25.023  7.404 25.070 /
\plot  7.404 25.070  7.355 25.118 /
\plot  7.355 25.118  7.305 25.169 /
\plot  7.305 25.169  7.254 25.224 /
\plot  7.254 25.224  7.201 25.279 /
\plot  7.201 25.279  7.148 25.339 /
\plot  7.148 25.339  7.093 25.400 /
\plot  7.093 25.400  7.038 25.459 /
\plot  7.038 25.459  6.983 25.521 /
\plot  6.983 25.521  6.930 25.580 /
\plot  6.930 25.580  6.879 25.637 /
\plot  6.879 25.637  6.833 25.688 /
\plot  6.833 25.688  6.790 25.737 /
\plot  6.790 25.737  6.756 25.777 /
\plot  6.756 25.777  6.727 25.809 /
\plot  6.727 25.809  6.703 25.836 /
\plot  6.703 25.836  6.687 25.855 /
\plot  6.687 25.855  6.676 25.866 /
\plot  6.676 25.866  6.670 25.874 /
\plot  6.670 25.874  6.668 25.876 /
\put{$x_n^\prime$} [lB] at  3.969 23.019
\put{$x_n$} [lB] at  3.969 21.273
\put{$\phi^{n}\left( h\left(  0\right)  \right)$} [lB] at  3.651 24.765
\put{$\phi\left( +\infty\right)$} [lB] at  1.111 22.066
\put{$\phi\left( - \infty\right)$} [lB] at 12.450 21.273
\put{$f(0)$} [lB] at  7.938 18.415
\put{$y_0^\prime$} [lB] at  7.779 20.003
\put{$y_0$} [lB] at  7.303 22.066
\put{$x_0$} [lB] at  9.366 21.114
\put{$h(0)$} [lB] at  9.684 23.654
\put{$\phi\left( h\left(  0\right)  \right)$} [lB] at  5.397 25.082
\linethickness=0pt
\putrectangle corners at  1.111 25.923 and 13.382 17.416
\endpicture}
\end{center}
\caption{}%
\label{telsx}%
\end{figure}

For this it suffices to find a geodesic $f^{\prime}$ in $\mathcal{O}$ such
that the image of $f^{\prime}$ intersects $Y_{0}$ (we then proceed to alter
$f^{\prime}$ and $K_{f^{\prime}}$ as we did with properties \textbf{(a)} and
\textbf{(b)} above). Let $Y_{-1}$ be the negatively curved polyhedron
constructed above. Since $Y_{-1}$ does not contain $1$-dimensional simplices,
the geodesic flow on $Y_{-1}$ is topologically mixing. Hence, the definition
\ref{mix} of topological mixing applied to the neighborhoods $\mathcal{O}%
_{-1}$ and $\mathcal{U}_{-1}^{\prime}$ (which are the images of the given
neighborhood $\mathcal{O}$ and the above chosen neighborhood $\mathcal{U}%
^{\prime}$ under the map (\ref{minusonesurj})) implies the existence of a
geodesic $f_{-1}^{\prime}\in\mathcal{O}_{-1}$ so that the image of
$f_{-1}^{\prime}$ intersects the subset of $Y_{-1}$ which corresponds to the
component $Y_{0}$. By property (\ref{minusonesurj}) above, choose a pre-image
$f^{\prime}\in\mathcal{O}$ of $f_{-1}^{\prime}.$ Then, by construction of
$Y_{-1},$ the image of $f^{\prime}$ intersects $Y_{0}$.

Thus, given arbitrary neighborhoods $\mathcal{O},\mathcal{U}\subseteq GY$ we
may choose geodesics $f^{\prime}\in\mathcal{O}$ and $h^{\prime}\in\mathcal{U}$
so that for some $t_{0}\in\mathbb{R},$ $f^{\prime}\left(  t\right)  $ and
$h^{\prime}\left(  t\right)  \in Y_{0}$ for all times $t\geq t_{0}.$ Let
$\widetilde{Y_{0}}$ be the universal cover of $Y_{0}.$
Since $\widetilde{Y_{0}}$ embeds isometrically in $\tilde{Y}$ we have, by
\cite[page 35]{C-D-P},  that
$\partial\widetilde{Y_{0}}$ injects in $\partial\tilde{Y} .$
Since $\Omega=GY$ (equivalently, $\Lambda\bigl(\pi
_{1}\left(  Y\right)  \bigr
)=\partial\tilde{Y})$ it follows that $\Lambda\bigl(\pi_{1}\left(
Y_{0}\right)  \bigr)=\partial\widetilde{Y_{0}}.$ As $\partial\widetilde{Y_{0}%
}$ is an infinite set, $\pi_{1}\left(  Y_{0}\right)  $ is non-elementary
which, together with property (\ref{5.1ofcoo}), implies that there exists a
hyperbolic isometry $\phi\in\pi_{1}\left(  Y_{0}\right)  \hookrightarrow
\pi_{1}\left(  Y\right)  $ such that for some lifts $f\in p^{-1}\left(
\mathcal{O}\right)  $ and $h\in p^{-1}\left(  \mathcal{U}\right)  $ of
$f^{\prime}$ and $h^{\prime},$ respectively, we have

\begin{itemize}
\item $\bigl(f\left(  +\infty\right)  ,h\left(  +\infty\right)  \bigr)=\left(
\phi\left(  +\infty\right)  ,\phi\left(  -\infty\right)  \right)  $,

\item  the image of the corresponding closed geodesic $c_{\phi}$ in $Y$ lies
entirely in $Y_{0}.$
\end{itemize}

We proceed now to show that for these choices of $f,h$ and $\phi$ the midpoint
of the geodesic segment $\left[  f\left(  0\right)  ,\phi^{n}\left(  h\left(
0\right)  \right)  \right]  $ lies, for all $n$ large enough, in the interior
of a $k$-simplex of $\tilde{Y}$ with $k\geq2.$ Then, proposition
\ref{alphazero} applies to assert the existence of $\xi_{n}$ in $\partial X$
such that
$\alpha\bigl(\xi_{n},f\left(  0\right)  ,\phi^{n}\left(  h\left(  0\right)
\right)  \bigr)=0$ completing the proof of \ref{mainncp}.

For each $n=0,1,2,\ldots,$ set $x_{n}$ to be the projection of $\phi
^{n}\left(  h\left(  0\right)  \right)  $ on the geodesic line $\left(
\phi\left(  +\infty\right)  ,\phi\left(  -\infty\right)  \right)  $ (by
notation, $\phi^{0}\left(  h\left(  0\right)  \right)  =h\left(  0\right)  $).
Similarly, let $y_{0}$ be the projection of $f\left(  0\right)  $ on the same
geodesic line. For the reader's convenience, we have gathered all the above
notation in figure 4. Let $\ell$ denote the translation length of the
hyperbolic isometry $\phi$ (i.e., $\ell$ equals the length of the closed
geodesic $\,c_{\phi}$ in $Y$ corresponding to $\phi$) and set
$\ell_{f}=d\left(  f\left(  0\right)  ,y_{0}\right)$ and
$\ell_{h}=d\left(  h\left(  0\right)  ,x_{0}\right).$
Observe that, since $\phi$ is an isometry, $d\bigl(\phi^{n}\left(  h\left(
0\right)  \right)  ,x_{n}\bigr)=\ell_{h}$ for all $n.$ Choose $N\in\mathbb{N}$
such that
\begin{equation}
n\ell>\ell_{f}+\ell_{h}+d\left(  x_{0},y_{0}\right)  \,\,\forall\,\,n\geq
N.\label{lproperty}%
\end{equation}
In addition, we may assume that
\begin{equation}
x_{n}\neq y_{0}\,\,\forall\,\,n\geq N.\label{pproperty}%
\end{equation}
The latter can be done because the sequence $\left\{  \phi^{n}\left(  h\left(
0\right)  \right)  \right\}  _{n\in\mathbb{N}}$ converges to $\phi\left(
+\infty\right)  $ and $d\bigl(\phi^{n}\left(  h\left(  0\right)  \right)
,x_{n}\bigr)=\ell_{h}$ for all $n.$

Fix any $n\geq N$ and consider
the piecewise geodesic segment
\[
\sigma:=\left[  f\left(  0\right)  ,y_{0}\right]  \cup\left[  y_{0}%
,x_{n}\right]  \cup\left[  x_{n},\phi^{n}\left(  h\left(  0\right)  \right)
\right].
\]
Consider the sub-path $\sigma^{\prime}$ of $\sigma$ which lies entirely in
$\widetilde{Y_{0}}$ and is of largest length with respect to this property.
Since $\left[  y_{0},x_{n}\right]  \subset\widetilde{Y_{0}},$ such a sub-path
exists and is of the form%
\[
\sigma^{\prime}=\left[  y_{0}^{\prime},y_{0}\right]  \cup\left[  y_{0}%
,x_{n}\right]  \cup\left[  x_{n},x_{n}^{\prime}\right]
\]
for some $y_{0}^{\prime}\in\left[  f\left(  0\right)  ,y_{0}\right]  $ and
$x_{n}^{\prime}\in\left[  x_{n},\phi^{n}\left(  h\left(  0\right)  \right)
\right]  .$

If $y_{0}^{\prime}=f\left(  0\right)  $ and $x_{n}%
^{\prime}=\phi^{n}\left(  h\left(  0\right)  \right)  ,$ it is clear that
$\sigma(=\sigma^{\prime})$  lies entirely in $\widetilde{Y_{0}}$ and  so does
the geodesic segment $\left[  f\left(  0\right)  ,\phi^{n}\left(  h\left(
0\right)  \right)  \right]  .$ In other words, $\left[  f\left(  0\right)
,\phi^{n}\left(  h\left(  0\right)  \right)  \right]  $ does not intersect
$\mathrm{Int\,}Y^{\left[  1\right]  },$ thus, its midpoint lies in the
interior of a $k$-simplex of $\tilde{Y}$ with $k\geq2$ as\pagebreak\ required.

If $y_{0}^{\prime}\neq f\left(  0\right)  $ and $x_{n}^{\prime}%
\neq\phi^{n}\left(  h\left(  0\right)  \right)  ,$ then $y_{0}^{\prime}$ and
$x_{n}^{\prime}$ are the $0$-faces of two distinct $1$-dimensional free simplices
attached to $\widetilde{Y_{0}}$ along $y_{0}^{\prime}$ and $x_{n}^{\prime}.$
The fact that $x_{n}^{\prime}$ and $y_{0}^{\prime}$ (as well as the
corresponding $1$-dimensional simplices) are distinct follows from
(\ref{pproperty}). It is easy to see that the piecewise geodesic
\[
\left[  f\left(  0\right)  ,y_{0}^{\prime}\right]  \cup\left[  y_{0}^{\prime
},x_{n}^{\prime}\right]  \cup\left[  x_{n}^{\prime},\phi^{n}\left(  h\left(
0\right)  \right)  \right]
\]
is, in fact, a geodesic, namely, is equal to the geodesic segment $\left[
f\left(  0\right)  ,\phi^{n}\left(  h\left(  0\right)  \right)  \right]  .$
Now property (\ref{lproperty}) asserts that the midpoint of $\left[  f\left(
0\right)  ,\phi^{n}\left(  h\left(  0\right)  \right)  \right]  $ lies in the
sub-segment $\left[  y_{0}^{\prime},x_{n}^{\prime}\right]  $ which is a subset
of $\widetilde{Y_{0}}.$

The cases $y_{0}^{\prime}=f\left(  0\right)
,x_{n}^{\prime}\neq\phi^{n}\left(  h\left(  0\right)  \right)  $ and
$y_{0}^{\prime}\neq f\left(  0\right)  ,x_{n}^{\prime}=\phi^{n}\left(
h\left(  0\right)  \right)  $ are treated in an identical way.\end{proof}

It is clear in the above proof that topological mixing in a polyhedron $Y$
which contains free $1$-simplices (i.e., $Y^{\left[  1\right]  }\neq\emptyset)$
is, in fact, deduced from topological mixing in the space $Y_{-1}$ which does
not contain free $1$-simplices. There is an alternative approach for this
deduction which uses (along with the space $Y_{-1}$) the topological mixing in
graphs as described in remark \ref{rfce}. This alternative approach
exhibits the important role of graphs in topological mixing of spaces of
higher dimension. Since the details are quite technical, we will only give a
description of this approach by presenting it in the frame of the following:

\begin{example}
Let $Y_{A},Y_{B}$ be two copies of a (triangulated) torus of genus 2
equip\textit{\-}ped with a metric of curvature $-1.$ Let $A\in Y_{A},$ $B\in
Y_{B}$ be $0$-simplices of $Y_{A}$ and $Y_{B}$ respectively. Let $Y$ be the
negatively curved polyherdron consisting of $Y_{A},Y_{B}$ and a single
$1$-simplex $\sigma$ (of length 1) whose $0$-faces are $A$ and $B,$ i.e.,
\[
Y=Y_{A}\vee_{A}\sigma\vee_{B}Y_{B}.
\]
In this example the space $Y_{-1}$ (defined in \ref{notation}) is the
one point union of $Y_{A}$ and $Y_{B}$ along $A$ and $B,$ i.e.,
\[
Y_{-1}=Y_{A}\vee_{A=B}Y_{B}.
\]
There are countably many geodesic segments in $Y_{A}$ starting and ending at
$A.$ Let $\Delta_{A}$ be the infinite graph consisting of

\begin{enumerate}
\item [-]a single vertex (denoted again by $A$),

\item[-] for each geodesic segment in $Y_{A}$ starting and ending at $A\in Y,$
there is one loop (of equal length) in $\Delta_{A}$ based at $A\in\Delta_{A}.$
\end{enumerate}\enlargethispage*{1\baselineskip}

Thus, $\Delta_{A}$ is an infinite graph (in fact, an infinite rose) whose
metric is determined by the geometry of $Y_{A}.$ Similarly, $\Delta_{B}$ is
defined. Let $\Delta$ be the disjoint union of $\Delta_{A}$ and $\Delta_{B}$
with one edge with endpoints $A$ and $B$ and of length 1 attached. Moreover,
set
\[
\Delta_{-1}=\Delta_{A} \vee_{A=B} \Delta_{B}.
\]
It is clear from the above construction that there is a bijection between
$G\Delta$ and the subset of $GY$ consisting of all geodesics $f\in GY$ such
that both $\operatorname{Im}\left(  f|_{\left[  0,+\infty\right)  }\right)  $
and $\operatorname{Im}\left(  f|_{\left(  -\infty,0\right]  }\right)  $
intersect $\sigma$ infinitely many times. By proposition \ref{alphazero} and
theorem \ref{main} the geodesic flow on $Y_{-1}$ is topologically mixing. It
follows that the lengths of all closed loops in $\Delta_{-1}$ generate a dense
subset of $\mathbb{R},$ otherwise, in a fashion similar to example \ref{ex}
and using the above-mentioned bijection, it would be possible to construct
neighborhoods in $GY_{-1}$ not satisfying definition \ref{mix}. Since closed
loops in $\Delta_{-1}$ generate a dense subset of $\mathbb{R},$ the same is
true for the closed loops in $\Delta,$ hence, topological mixing holds for
$\Delta$ (cf. remark \ref{rfce}). On the other\pagebreak\ hand, it can be shown that
topological mixing on $\Delta$ implies (in fact, is equivalent with)
topological mixing on $Y.$ 

In conclusion, topological mixing in $Y$
can proved through graphs once we know that it holds for $Y_{-1}.$
This alternative approach described in this example can be defined and proved
in detail for an arbitrary negatively curved polyhedron $Y.$
\end{example}

\begin{remark}
Let $Y$ be a negatively curved polyhedron; then $Y$ is a graph precisely when
$\partial\tilde{Y}$ is totally disconnected. Thus, theorem \ref{mainncp} says
that the geodesic flow on a negatively curved polyhedron is topologically
mixing if the boundary of its universal cover is not totally disconnected. It
is plausible to expect that this is the case in any geodesically complete
$CAT\left(  -1\right)  $-space, namely, that condition (1) in theorem
\ref{main} can be replaced by the assumption that $\partial\tilde{Y}$ is not
totally disconnected.
\end{remark}

\subsection{Ideal polyhedra}

We apply theorem \ref{main} to a non-compact class of spaces, namely, to
$n$-\textit{dimen\allowbreak sional complete ideal polyhedra}. Important
examples of ideal polyhedra have appeared in Thurston's work (see \cite{Thu},
\cite[Sec.\ 10.3]{Rat}), where $3$-mani\textit{\-}folds, which are complements of
links and knots in $\mathbb{S}^{3}$, are constructed by gluing together
finitely many ideal tetrahedra. In consequence, these finite volume
$3$-manifolds are equipped by a complete hyperbolic structure. Moreover, the
$2$-skeleton of these $3$-manifolds are examples of $2$-dimensional ideal
polyhedra. 

\begin{definition}
An $n$-dimensional ideal polyhedron is a locally finite union of ideal
hyperbolic $n$-polytopes glued together isometrically along their $\left(
n-1\right)  $-faces with at least two germs of polytopes along each $\left(
n-1\right)  $-face. The distance function is defined exactly as 
described in section \ref{negpolsec} for negatively curved polyhedra, i.e., the
distance $d\left(  x,y\right)  $ from $x$ to $y$ is defined to be the lower
bound of the lengths of broken geodesics from $x$ to $y.$ With the induced
metric, an $n$-dimensional ideal polyhedron is required to be complete. Since
it is locally compact, it is proper and geodesic. Moreover, it is required to
have curvature less than or equal to $-1.$
\end{definition}

We note here that in the case $n=2$ the curvature condition in the
above definition can be proved, hence, is redundant (see \cite[Prop.\ 1]{Ch}).
An ideal polyhedron $Y$ is called \textit{finite} if finitely many polytopes
are glued together to form $Y.$

Certain properties of this class of spaces, including transitivity of the
geodesic flow, have been studied in \cite{C-P}, \cite{C-T}, \cite{Ch} and
\cite{C-T1}.

If $Y$ is an ideal polyhedron of dimension $n$, then $Y$ is naturally a proper
geodesic metric space.

The universal covering $\widetilde{Y}$ of $Y$ is a complete ideal polyhedron
of dimension $n$ satisfying $CAT\left(  -1\right)  $ inequality (see
\cite[Cor. 2.11]{Pau}). If $Y$ is a finite polyhedron, the non-wandering set
$\Omega$ of the geodesic flow on $Y$ is equal to $GY$ (see Cor. 10 in
\cite{C-T}) and $\pi_{1}\left(  Y\right)  $ is a non-elementary group of
isometries acting properly discontinuously on $\widetilde{Y}$ (see Cor. 12 in
\cite{C-T}). Moreover, proposition \ref{alphazero} applies verbatim to ideal
polyhedra. Hence, we obtain the following application of theorem \ref{main}.

\begin{corollary}
Let $Y$ be an $n$-dimensional finite ideal polyhedron. Then the geodesic flow
on $Y$ is topologically mixing.
\end{corollary}

\section*{Acknowledgments}
We would like to thank the referee for very
helpful suggestions and observations which initiated significant improvements.
One of them was the removal of a restrictive hypothesis in the main theorem
and the enhancement of section \ref{application}\pagebreak.

\bibliographystyle{amsplain}

\begin{thebibliography}{99}

\bibitem{Ano}D.V. Anosov, \textit{Geodesic flows on closed Riemannian
manifolds with negative curvature}, Proc. Steklov Instit. Math. 90, Amer.
Math. Soc., Providence, Rhode Island, 1969. \MR{39:3527}

\bibitem{Bal}W. Ballmann, \textit{Lectures on Spaces of Non-positive
Curvature}, Birkh\"{a}user, 1995. \MR{97a:53053}

\bibitem{B-G-H}W. Ballman, E. Ghys, A. Haefliger, P. de la Harpe, E. Salem, R.
Strebel and M. Troyanov, \textit{Sur les groups hyperboliques d'apr\'{e}s
Gromov} (Seminaire de Berne), \'{e}dit\'{e} par E. Ghys et P. de la Harpe, (a
paraitre chez Birkh\"{a}user), 1990. 

\bibitem{Bou}M. Bourdon, \textit{Structure conforme au bord et flot
g\'eod\'esique d'un }CAT$\left(  -1\right)  $-e\textit{space}, Enseign. Math.
41, (1995), 63-102. \MR{96f:58120}

\bibitem{Bow1}B. H. Bowditch, \textit{Boundaries of geometrically finite
groups, }Math. Z. 230 (1999) no. 3, 509-527. \MR{2000b:20049}

\bibitem{Bow2}B. H. Bowditch, \textit{Connectedness properties of limit
sets}, Trans. Amer. Math. Soc. 351 (1999) 3673-3686. \MR{2000d:20056}

\bibitem{Bri}M. Bridson, \textit{Geodesics and curvature in metric simplicial
complexes} in Group Theory from a Geometrical Viewpoint, (ICTP, Trieste,
Italy, March 26-April 6, 1990), E.Ghys and A.Haefliger, eds., (1991).
\MR{94c:57040}

\bibitem{Cham}Ch. Champetier, \textit{Petite simplification dans les groupes
hyperboliques}, Ann. Fac. Sci. Toulouse, VI. Ser., Math. 3, No.2, 161-221,
(1994). \MR{95e:20050}

\bibitem{C-P}C. Charitos, A. Papadopoulos, \textit{The geometry of ideal
polyhedra}, to appear in Glasgow Journal of Mathematics.

\bibitem{C-T}C. Charitos, G. Tsapogas, \textit{Geodesic flow on ideal
polyhedra}, Canad. J. Math. 49 (4) 1997, pp. 696-707. \MR{98e:58133}

\bibitem{Ch}C. Charitos, \emph{Closed geodesics in ideal polyhedra of
dimension} 2, Rocky Mountain Journal of Mathematics, Vol. 26, no. 2 (1996), pp.
507-521. \MR{97d:57014}

\bibitem{C-T1}C. Charitos, G. Tsapogas, \textit{Complexity of Geodesics on
2-dimensional ideal polyhedra and isotopies}, Math. Proc. Camb. Phil. Soc.,
Vol. 121 (1997), pp. 343-358. \MR{98b:57006}

\bibitem{C-T2}C. Charitos, G. Tsapogas, \textit{Approximation of recurrence in
negatively curved metric spaces,} Pacific J. Math. 195 (2000), 67--79.


\bibitem{Coo1}M. Coornaert, \textit{Sur les groupes proprement discontinus
d'isom\'{e}tries des espaces hyperboliques au sens de Gromov}, Th\`{e}se
U.L.P., Publication de l'IRMA. \MR{92i:57032}

\bibitem{C-D-P}M. Coornaert, T. Delzant, A. Papadopoulos, \textit{G\'{e}%
om\'{e}trie et th\'{e}orie des groupes}, Lecture Notes in Mathematics,
vol.1441, Springer-Verlag, (1990). \MR{92f:57003}

\bibitem{Ebe2}P. Eberlein, \textit{Geodesic flows on negatively curved
manifolds II}, Trans. Amer. Math. Soc. 178 (1973), pp. 57-82. \MR{47:2636}

\bibitem{Gro1}M. Gromov, \textit{Hyperbolic groups}, in Essays in Group
Theory, MSRI Publ. 8, Springer, 1987, pp. 75-263. \MR{89e:20070}

\bibitem{H-P}S. Hersonsky, F. Paulin, \textit{On the rigidity of discrete
isometry groups of negatively curved spaces}, Comm. Math. Helv. 72 (1997),
pp. 349-388. \MR{98h:58105}

\bibitem{Kai}V.A. Kaimanovich, \textit{Ergodicity of harmonic invariant
measures for the geodesic flow on hyperbolic spaces}, J. Reine Angew. Math.
445 (1994), pp. 57-103. \MR{95g:58130}

\bibitem{Mou}G. Moussong, \textit{Hyperbolic Coxeter groups}, Doctoral
Dissertation, Ohio State University, 1988.

\bibitem{Pau}F. Paulin, \textit{Constructions of hyperbolic groups via
hyperbolization of polyhedra}, in Group Theory from a Geometrical
Viewpoint, (ICTP, Trieste, Italy, March 26-April 6, 1990), E. Ghys and A.
Haefliger, eds., (1991). \MR{93d:57005}

\bibitem{Rat}J. Ratcliffe, \textit{Foundations of hyperbolic geometry}, GTM,
Springer-Verlag, 1994. \MR{95j:57011}

\bibitem{Thu}W.P. Thurston, \textit{The Geometry and Topology of
Three-manifolds}, Lecture notes, Princeton University, Princeton, NJ (1979).
\enlargethispage*{2\baselineskip}
\end{thebibliography}

\end{document}